# Explicit maximal totally real embeddings

BY

**NEFTON PALI** and **BRUNO SALVY**

**Abstract**

This article is the continuation of the first named author work "On maximal totally real embeddings". For real analytic compact manifolds equipped with a covariant derivative operator acting on the real analytic sections of its tangent bundle, a construction of canonical maximal totally real embeddings is known from previous works by Guillemin-Stenzel, Lempert, Lempert-Szöke, Szöke and Bielawski. The construction is based on the use of Jacobi fields, which are far from being explicit. As a consequence, the form of the corresponding complex structure has been a mystery since the very beginning. A quite simple recursive expression for such complex structures has been provided in the above cited first named author work. In our series of papers we always focus on the torsion free case. In the present paper we provide a fiberwise Taylor expansion of the canonical complex structure which is expressed in terms of symmetrization of curvature monomials and a rather simple and explicit expression of the coefficients of the expansion. Our main argument applies to far more generals settings that can be useful for the study of open questions in the theory of the embeddings in consideration. In this paper we provide also evidence for some remarkable canonical vanishing of some of the integrability equations in general settings.

## 1 Introduction and statement of the main result

**Notations and definitions for the statement of the main theorem.**

We remind first some notations and definitions from the work [Pali]. We remind first that over a smooth manifold $M$, the vector bundle $T_{T_M|M} \simeq T_M \oplus T_M$, is a complex one with the canonical complex structure $J^{\mathrm{can}} : (u, v) \longmapsto (-v, u)$ acting on the fibers. Any almost complex structure which is a continuous extension of $J^{\mathrm{can}}$ in a neighborhood of $M$ inside $T_M$ makes $M$ an (almost complex) maximally totally real sub-manifold of $T_M$. We denote by

$$T \in C^{\infty}(T_M, \pi^* T_M^* \otimes_{\mathbb{R}} T_{T_M}),$$

the canonical section which at the point $\eta \in T_M$ takes the value $T_\eta : E_p \longrightarrow T_{E_p, \eta}$, which is the canonical isomorphism map.

**Definition 1.1.** *Let $M$ be a smooth manifold. An $M$-totally real almost complex structure over an open neighborhood $U \subseteq T_M$ of the zero section is a couple $(\alpha, B)$ with $\alpha \in C^{\infty}(U, \pi^* T_M^* \otimes_{\mathbb{R}} T_{T_M})$ and $B \in C^{\infty}(U, \pi^* \mathrm{GL}(T_M))$ such that $d\pi \cdot \alpha = \mathbb{I}_{\pi^* T_M}$ over $U$ and such that $\alpha_{0_p} = d_p 0_M$, $B_{0_p} = \mathbb{I}_{T_M, p}$, for all $p \in M$. The almost complex structure $J_A$, with $A = \alpha + iTB$ associated to $(\alpha, B)$ is the one which satisfies $T^{0,1}_{T_M, J_A, \eta} = A_\eta(\mathbb{C}T_{M, \pi(\eta)}) \subset \mathbb{C}T_{T_M, \eta}$ for all $\eta \in U \subseteq T_M$.*

Every almost complex continuous extension of the canonical complex structure $J^{\mathrm{can}}$ of $T_{T_M|M}$ over a neighborhood of $M$ inside $T_M$ writes as the almost complex structure associated to an $M$-totally real almost complex structure defined over a sufficiently small neighborhood of $M$.

**Notation 1.2.** *We will denote by* Circ *the circular operator*

$$(\mathrm{Circ}\,\theta)(v_1, v_2, v_3, \bullet) = \theta(v_1, v_2, v_3, \bullet) + \theta(v_2, v_3, v_1, \bullet) + \theta(v_3, v_1, v_2, \bullet),$$

*acting on the first three entries of any $q$-tensor $\theta$, with $q \geqslant 3$. We define also the permutation operation $\theta_2(v_1, v_2, \bullet) := \theta(v_2, v_1, \bullet)$.*

*For any complex covariant derivative $\nabla$ operator acting on the smooth sections of $\mathbb{C}T_M$ we define the operator*

$$d_1^{\nabla} \colon C^{\infty}(M, T_M^{*,\otimes k} \otimes_{\mathbb{R}} \mathbb{C}T_M) \longrightarrow C^{\infty}\big(M, \Lambda^2 T_M^* \otimes_{\mathbb{R}} T_M^{*,\otimes(k-1)} \otimes_{\mathbb{R}} \mathbb{C}T_M\big),$$

*with $k \geqslant 1$ as follows*

$$d_1^{\nabla} A(\xi_1, \xi_2, \mu) := \nabla_{\xi_1} A(\xi_2, \mu) - \nabla_{\xi_2} A(\xi_1, \mu),$$

*with $\xi_1$, $\xi_2 \in T_M$ and with $\mu \in T_M^{\oplus(k-1)}$. For notation simplicity we will denote by $d_1^{\nabla,p} := (d_1^{\nabla})^p$. Moreover for any $A \in C^{\infty}(M, T_M^{*,\otimes(k+1)} \otimes_{\mathbb{R}} \mathbb{C}T_M)$ and $B \in C^{\infty}(M, T_M^{*,\otimes(l+1)} \otimes_{\mathbb{R}} \mathbb{C}T_M)$ we define the exterior product*

$$A \wedge_1 B \in C^{\infty}\big(M, \Lambda^2 T_M^* \otimes_{\mathbb{R}} T_M^{*,\otimes(k+l-1)} \otimes_{\mathbb{R}} \mathbb{C}T_M\big),$$

*as*

$$(A \wedge_1 B)(\xi_1, \xi_2, \eta, \mu) := A(\xi_1, B(\xi_2, \eta), \mu) - A(\xi_2, B(\xi_1, \eta), \mu),$$

*with $\xi_1$, $\xi_2 \in T_M$, $\eta \in T_M^{\oplus l}$ and $\mu \in T_M^{\oplus(k-1)}$. We denote by $\mathrm{Sym}_{r_1,\dots,r_s}$ the symmetrizing operator (without normalizing coefficient!) acting on the entries $r_1$, ..., $r_s$ of a multi-linear form. We use in this paper the common convention that a sum and a product running over an empty set is equal respectively to 0 and 1.*

We observe that the argument given in the proof of corollary 2 in [Pali] applies without modifications in the case of a general torsion free covariant derivative operator acting on the smooth sections of the tangent bundle. For the readers convenience we will provide in the appendix, the argument in this more general setting. We have therefore the following statement.

**Corollary 1.3.** *Let $M$ be a smooth manifold, let $\nabla$ be a torsion free covariant derivative operator acting on the smooth sections of $T_M$ and let $J \equiv J_A$ be an $M$-totally real almost complex structure over an open neighborhood $U$ of $M$ inside $T_M$, with connected fibers, which is real analytic along the fibers of $U \subset T_M$. Then the structure $J$ is integrable over $U$ and for any $\eta \in U$, the complex curve $\psi_\eta : t + is \longmapsto s\Phi_t^{\nabla}(\eta)$, defined in a neighborhood of $0 \in \mathbb{C}$, is $J$-holomorphic if and only if the fiberwise Taylor expansion at the origin*

$$T_\eta^{-1}(H^{\nabla} - \overline{A})_\eta \cdot \xi = i\xi + \sum_{k \geqslant 1} S_k(\xi, \eta^k),$$

*for any $\eta \in T_M$ in a sufficiently small neighborhood of the zero section and for any $\xi \in T_{M,\pi(\eta)}$ with $S_k \in C^{\infty}(M, T_M^* \otimes_{\mathbb{R}} S^k T_M^* \otimes_{\mathbb{R}} \mathbb{C}T_M)$ and with $\eta^k := \eta^{\times k} \in T_{M,\pi(\eta)}^{\oplus k}$, satisfies: $S_1 = 0$,*

$$S_k = \frac{i}{(k+1)!k!} \mathrm{Sym}_{2,\dots,k+1} \Theta_k^{\nabla},$$

*for all $k \geqslant 2$, with $\Theta_2^{\nabla} := 2R^{\nabla}$, with*

$$\begin{aligned} \Theta_k^{\nabla} &:= -2i\,(i d_1^{\nabla})^{k-3} (\nabla R^{\nabla})_2 \\ &+ \sum_{r=3}^{k-1} (r+1)! \sum_{p=2}^{r-1} (i d_1^{\nabla})^{k-1-r} (p S_p \wedge_1 S_{r-p+1}), \end{aligned}$$

*for all $k \geqslant 3$ and $\mathrm{Circ}\,\mathrm{Sym}_{3,\dots,k+1} \Theta_k^{\nabla} = 0$, for all $k \geqslant 4$.*

We have shown in [Pali] that the equation $\operatorname{Circ}\operatorname{Sym}_{3,\ldots,k+1}\Theta_k^{\nabla}=0$ is satisfied for $k=4$ even in the more general case of a torsion free complex covariant derivative operator acting on the smooth sections of the complexified tangent bundle. We wish to point out that in this more general case it is not possible to have a statement similar to corollary 1.3, simply because there are no geodesics in $M$ associated to the complex covariant derivative operator. However this complex case is rather important to the applications to micro-local analysis since quite often the operators there are expressed in terms of complex covariant derivative operators. In this paper we build a general formalism which allows us to give a more compact proof of the vanishing of the tensor $\operatorname{Circ}\operatorname{Sym}_{3,\ldots,k+1}\Theta_k^{\nabla}$ in the case $k=4$, and which allows us to provide also a proof in the case $k=5$. We have checked by using computer programing the above vanishing also for the cases $k=6,7$. It is quite likely that the vanishing is always satisfied in the torsion free complex covariant derivative case. By corollary 1.3, this vanishing is certainly true in the cases considered by Guillemin-Stenzel [Gu-St], Lempert [Lem], Lempert-Szöke [Le-Sz1, Le-Sz2], Szöke [Szo] and Bielawski [Bie]. We feel however that a proof independent of their work should provide a more general statement.

We provide now some more notations useful for the statement of our main theorem.

**Notation 1.4.** *We denote by* $A\vDash_0 p$ *any element* $A\equiv(a_1,\ldots,a_N)\in\mathbb{Z}_{\geqslant 0}^N$ *with* $N\in\mathbb{Z}_{>0}$ *arbitrary such that* $|A|:=\sum_{j=1}^N a_j=p$*. We set* $l_A:=N$*. We denote also*

$$|A|_j := \sum_{s=1}^{j} a_s\,,$$

$$\|A\|_j := |A|_j+j\,,$$

*for any* $1\leqslant j\leqslant N$*, and we set* $\|A\|:=\|A\|_N+N$*. For any* $l=1,\ldots,N$*, we define* $A_l^-:=(a_1,\ldots,a_l)$ $A_l^+:=(a_{l+1},\ldots,a_N)$ *and*

$$A'_{l,h} := (A_h^-)_l^+=(a_{l+1},\ldots,a_h)\,,$$

*any* $l,h\in\mathbb{Z}_{\geqslant 0}$ *with* $l\leqslant h\leqslant N$*. In all this paper, for notations simplicity, we will use the identification* $(v_1,\ldots,v_p)\equiv(1,\ldots,p)$*, for a given element* $(v_1,\ldots,v_p)\in T_M^{\oplus p}$*.*

*Let* $R^{(a)}:=\nabla^a R^{\nabla}$*. We denote by* $R^{(a)}(1,\ldots,a+1,\bullet,a+2)$ *the endomorphism*

$$v \longmapsto R^{(a)}(1,\ldots,a+1,v,a+2)\,.$$

*Given a family of endomorphisms* $(T_j)_{j=1}^p$ *we denote by*

$$\prod_{j=1}^{p} T_j := T_1\circ\cdots\circ T_p\,.$$

*We define also*

$$R^A(1,\ldots,\|A\|+1) := \left[\prod_{j=1}^{l_A} R^{(a_j)}\big(\|A\|_{j-1}+2\,,\ldots,\|A\|_j+1,\bullet,\|A\|-j+2\big)\right]\cdot 1\,,$$

*i.e. the "monomial"* $R^A$ *writes as*

$$\begin{aligned}&R^A(1,\ldots,\|A\|+1)\\ =\;& R^{(a_1)}\big(2,\ldots,a_1+2,R^{(a_2)}\big(a_1+3,\ldots,a_1+a_2+3,R^{(a_3)}\big(a_1+a_2+4,\ldots,a_1+a_2+a_3+4,R^{(a_4)}\big(\ldots\\ &\ldots,R^{(a_{l_A})}\big(\|A\|_{l_A-1}+2,\ldots,\|A\|_{l_A}+1,1,\|A\|_{l_A}+2\big),\|A\|_{l_A}+3\big)\,,\ldots\big),\|A\|+1\big)\,.\end{aligned}$$

*We denote by* $\lambda \vDash l$ *any* $\lambda \in \mathbb{Z}^p_{\geqslant 1}$ *such that* $|\lambda| = l$ *and we define for any* $A \in \mathbb{Z}^l_{\geqslant 0}$ *the coefficient*

$$C_A \;:=\; -(-i)^{\|A\|}\,(\|A\|+1)! \sum_{0 \leqslant H \leqslant A} \frac{(-1)^{|H|}\,C(H)}{H!\,(A-H)!}\,, \quad \textit{with} \quad l_H = l\,,$$

$$C(H) \;:=\; \sum_{\lambda \vDash l} (-1)^{l_\lambda} \|H^-_{\lambda_1}\| \prod_{j=1}^{l_\lambda} \prod_{s=|\lambda|_{j-1}}^{|\lambda|_j - 1} \frac{1}{\|H'_{s,|\lambda|_j}\|(\|H'_{s,|\lambda|_j}\|+1)}\,.$$

With theese notations we can state our main theorem.

**Theorem 1.5.** (***Main theorem: Explicit maximal totally real embeddings***)
*Let* $M$ *be a smooth manifold, let* $\nabla$ *be a torsion free covariant derivative operator acting on the smooth sections of* $T_M$ *and let* $J \equiv J_A$ *be an* $M$*-totally real almost complex structure over an open neighborhood* $U$ *of* $M$ *inside* $T_M$*, with connected fibers, which is real analytic along the fibers of* $U \subset T_M$*. Then* $J$ *is integrable over* $U$ *and for any* $\eta \in U$*, the complex curve* $\psi_\eta : t + is \longmapsto s\,\Phi^\nabla_t(\eta)$*, defined in a neighborhood of* $0 \in \mathbb{C}$*, is* $J$*-holomorphic if and only if the fiberwise Taylor expansion at the origin*

$$T^{-1}_\eta (H^\nabla - \overline{A})_\eta \cdot \xi \;=\; i\xi + \sum_{k \geqslant 1} S_k(\xi, \eta^k)\,,$$

*for any* $\eta \in T_M$ *in a sufficiently small neighborhood of the zero section and for any* $\xi \in T_{M,\pi(\eta)}$ *with* $S_k \in C^\infty(M, T^*_M \otimes_{\mathbb{R}} S^k T^*_M \otimes_{\mathbb{R}} \mathbb{C}T_M)$ *and with* $\eta^k := \eta^{\times k} \in T^{\oplus k}_{M,\pi(\eta)}$*, satisfies:* $S_1 = 0$*,*

$$S_k \;=\; \frac{i}{(k+1)!\,k!} \sum_{\|D\|=k} C_D\, \mathrm{Sym}_{2,\ldots,k+1} R^D\,, \quad \textit{with} \quad D \geqslant 0,$$

*for all* $k \geqslant 2$ *and* $\mathrm{Circ}\,\mathrm{Sym}_{3,\ldots,k+1}\,\Theta^\nabla_k = 0$*, for all* $k \geqslant 6$*, where*

$$\Theta^\nabla_k \;:=\; -2i\,(i d^\nabla_1)^{k-3}(\nabla R^\nabla)_2 + \sum_{r=3}^{k-1} (r+1)! \sum_{p=2}^{r-1} (i d^\nabla_1)^{k-1-r}(p S_p \wedge_1 S_{r-p+1})\,.$$

The main theorem is a direct consequence of corollary 1.3, combined with the following general result (as well as the vanishing for $k = 4, 5$ that we will provide below).

**Theorem 1.6.** *Let* $M$ *be a smooth manifold and let* $\nabla$ *be a torsion free complex covariant derivative operator acting on the smooth sections of the complexified tangent bundle* $\mathbb{C}T_M$*. Then*

$$S_k \;=\; \frac{i}{(k+1)!\,k!} \mathrm{Sym}_{2,\ldots,k+1}\,\Theta^\nabla_k\,,$$

*for all* $k \geqslant 2$*, with* $\Theta^\nabla_2 := 2R^\nabla$*, with*

$$\Theta^\nabla_k \;:=\; -2i\,(i d^\nabla_1)^{k-3}(\nabla R^\nabla)_2 + \sum_{r=3}^{k-1} (r+1)! \sum_{p=2}^{r-1} (i d^\nabla_1)^{k-1-r}(p S_p \wedge_1 S_{r-p+1})\,,$$

*for all* $k \geqslant 3$ *if and only if*

$$S_k \;=\; \frac{i}{(k+1)!\,k!} \sum_{\|D\|=k} C_D\, \mathrm{Sym}_{2,\ldots,k+1} R^D\,, \quad \textit{with} \quad D \geqslant 0,$$

*for all* $k \geqslant 2$.

In the auxiliary preprint [Pal-Sal], we provide a Maple program verifying the main equivalent statement in the proof of the theorem 1.6. We warmly thank François Guenard for providing an alternative verification of the same statement using the software "Mathematica", that we also include there.

# 2 Proof of theorem 1.6

In all the paper, with the exception of the appendix, we assume $\nabla$ be a torsion free complex covariant derivative operator acting on the smooth sections of the complexified tangent bundle $\mathbb{C}T_M$.

## 2.1 Expliciting the powers of the 1-differential. Part I

**Notation 2.1.** *For any* $A \in T_M^{*,\otimes p} \otimes \mathrm{End}_{\mathbb{C}}(\mathbb{C}T_M)$ *and for any* $\theta \in T_M^{*,\otimes q} \otimes \mathbb{C}T_M$*, the product operations of tensors* $A \cdot \theta, A \neg \theta \in T_M^{*,\otimes (p+q)} \otimes \mathbb{C}T_M$ *are defined by*

$$(A \cdot \theta)(u_1, ..., u_p, v_1, ..., v_q) := A(u_1, ..., u_p) \cdot \theta(v_1, ..., v_q),$$

$$(A \neg \theta)(u_1, ..., u_p, v_1, ..., v_q) := \sum_{j=1}^{q} \theta(v_1, ..., A(u_1, ..., u_p) \cdot v_j, ..., v_q).$$

*We will denote for notation simplicity* $R^{\nabla}.\theta := R^{\nabla} \cdot \theta - R^{\nabla} \neg \theta$.

We denote by $\mathrm{Alt}_2$, the alternating operator (without normalizing coefficient!) acting on the first tow entries of any tensor. We remind the following well known fact (see [Pali]).

**Lemma 2.2.** *For any complex covariant derivative operator* $\nabla$ *acting on the smooth sections of the complexified tangent bundle* $\mathbb{C}T_M$ *and for any tensor* $\theta \in C^{\infty}(M, T_M^{*,\otimes k} \otimes_{\mathbb{R}} \mathbb{C}T_M)$*, holds the commutation identity*

$$\mathrm{Alt}_2 \nabla^2 \theta = R^{\nabla}.\theta. \tag{2.1}$$

We observe now other two elementary lemmas. Their proof is left to the reader.

**Lemma 2.3.** *For any* $p$*-tensor* $\theta$ *hold the formula*

$$(d_1^{\nabla})^k \theta (1, ..., k+p) = \sum_{\sigma \in S_{k+1}^*} \varepsilon_{\sigma} \nabla^k \theta (\sigma_1, ..., \sigma_{k+1}, k+2, ..., k+p),$$

*where* $S_{k+1}^*$ *is the set of permutations* $\sigma$ *of the set* $\{1, ..., k+1\}$ *such that if* $1 \leqslant j < k < l \leqslant k+1$ *then* $\sigma_j < \sigma_k$ *and* $\sigma_j < \sigma_l$.

**Lemma 2.4.** *For any* $p$*-tensor* $\theta$ *hold the formula*

$$\begin{aligned}
&(d_1^{\nabla})^k \theta (1, ..., k+p)\\
=& \sum_{\substack{j_{r-1} < l_r < j_r \leqslant k+1\\ j_0 := 0, r \geqslant 1}} (-1)^{\sum_{r \geqslant 1} (j_r - l_r)} \nabla^k \theta (1, ..., \hat{l}_1, ..., j_1, l_1, j_1 + 1, ...,\\
&, ..., \hat{l}_2, ..., j_2, l_2, j_2 + 1, ..., \hat{l}_r, ..., j_r, l_r, j_r + 1, ..., k+2, ..., k+p).
\end{aligned}$$

We remind that the notation $a, ..., b$ for integers $a < b$ denotes the increasing by one sequence from the left to the right hand side. We ignore this notation when $a > b$. We start by showing the following fundamental proposition.

**Proposition 2.5.** *Let $\nabla$ be a torsion free complex covariant derivative operator acting on the smooth sections of $\mathbb{C}T_M$ with curvature operator $R^\nabla(\cdot,\cdot)\cdot \equiv R^\nabla(\cdot,\cdot,\cdot)$. Then for all integers $k \geqslant 1$ hold the identities*

$$\operatorname{Sym}_{2,...,k+4}\Big[(d_1^\nabla)^k\,(\nabla R^\nabla)_2\Big] = \operatorname{Sym}_{2,...,k+4}\,\Phi_k\,, \tag{2.2}$$

$$\operatorname{Sym}_{3,...,k+4}\Big[(d_1^\nabla)^k\,(\nabla R^\nabla)_2\Big] = \operatorname{Sym}_{3,...,k+4}(T_k + Q_k + V_k)\,, \tag{2.3}$$

*with*

$$\begin{aligned}
\Phi_k(1,...,k+4) &:= -(-1)^k\,\nabla^{k+1}R^\nabla(2,...,k+3,1,k+4)\,,\\
T_k(1,...,k+4) &:= (-1)^k\,\nabla^{k+1}R^\nabla(3,...,k+3,1,2,k+4)\,,\\
Q_k(1,...,k+4) &:= -(-1)^k\sum_{j=2}^{k+1}\nabla^{j-2}\big(R^\nabla.\nabla^{k+1-j}R^\nabla\big)(3,...,j,2,j+1,...,k+3,1,k+4)\,,\\
V_k(1,...,k+4) &:= (-1)^k\sum_{j=2}^{k+1}\nabla^{j-2}\big(R^\nabla.\nabla^{k+1-j}R^\nabla\big)(3,...,j,1,j+1,...,k+3,2,k+4)\,.
\end{aligned}$$

**Proof.** We notice first that the symmetrization of a tensor with two alternating entries vanishes. Applying this fact to the covariant derivatives of $R^\nabla$ we infer the identities

$$\begin{aligned}
\operatorname{Sym}_{2,...,k+4}\Big[(d_1^\nabla)^k\,(\nabla R^\nabla)_2\Big] &= \operatorname{Sym}_{2,...,k+4}\,\varphi_k\,,\\
\operatorname{Sym}_{3,...,k+4}\Big[(d_1^\nabla)^k\,(\nabla R^\nabla)_2\Big] &= \operatorname{Sym}_{3,...,k+4}(\varphi_k + \psi_k)\,,
\end{aligned}$$

with

$$\begin{aligned}
\varphi_k(1,...,k+4) &:= (-1)^k\nabla^k(\nabla R^\nabla)_2\,(2,...,k+1,1,k+2,k+3,k+4)\,,\\
\psi_k(1,...,k+4) &:= -(-1)^k\nabla^k(\nabla R^\nabla)_2\,(1,3,...,k+1,2,k+2,k+3,k+4)\,,
\end{aligned}$$

which rewrite as

$$\begin{aligned}
\varphi_k(1,...,k+4) &= (-1)^k\,\nabla^{k+1}R^\nabla\,(2,...,k+2,1,k+3,k+4)\\
&= -(-1)^k\,\nabla^{k+1}R^\nabla\,(2,...,k+3,1,k+4)\,,\\
\psi_k(1,...,k+4) &= -(-1)^k\,\nabla^{k+1}R^\nabla\,(1,3,...,k+2,2,k+3,k+4)\\
&= (-1)^k\,\nabla^{k+1}R^\nabla\,(1,3,...,k+3,2,k+4)\,,
\end{aligned}$$

thanks to the alternating property of $R^{\nabla}$. We deduce in particular the identity (2.2). Using the identity (2.1) we infer

$$\begin{aligned}
\varphi_k(1,...,k+4) &= -(-1)^k\, \nabla^{k+1}R^{\nabla}\,(3,2,4,...,k+3,1,k+4)\\
&- (-1)^k\big(R^{\nabla}.\nabla^{k-1}R^{\nabla}\big)(2,...,k+3,1,k+4)\,,\\
\psi_k(1,...,k+4) &= (-1)^k\, \nabla^{k+1}R^{\nabla}\,(3,1,4,...,k+3,2,k+4)\\
&+ (-1)^k\big(R^{\nabla}.\nabla^{k-1}R^{\nabla}\big)(1,3,...,k+3,2,k+4)\,.
\end{aligned}$$

We show now for any integer $p$, with $2\leqslant p\leqslant k+1$, the identities

$$\begin{aligned}
\varphi_k(1,...,k+4) &= -(-1)^k\nabla^{k+1}R^{\nabla}(3,...,p+1,2,p+2,...,k+3,1,k+4)\\
&- (-1)^k\sum_{j=2}^{p}\nabla^{j-2}\big(R^{\nabla}.\nabla^{k+1-j}R^{\nabla}\big)(3,...,j,2,j+1,...,k+3,1,k+4)\,,\\
\psi_k(1,...,k+4) &= (-1)^k\, \nabla^{k+1}R^{\nabla}\,(3,...,p+1,1,p+2,...,k+3,2,k+4)\\
&+ (-1)^k\sum_{j=2}^{p}\nabla^{j-2}\big(R^{\nabla}.\nabla^{k+1-j}R^{\nabla}\big)(3,...,j,1,j+1,...,k+3,2,k+4)\,.
\end{aligned}$$

We show them by finite induction on $p$. We assume them true for $p<k+1$. Applying the covariant derivative $\nabla^{p-1}_{3,...,p+1}$ to the identity (2.1), with $\theta:=\nabla^{k+2-p}R^{\nabla}$, we infer the conclusion of the induction. If we set $p:=k+1$ in the previous identites we obtain

$$\begin{aligned}
(\varphi_k+\psi_k)(1,...,k+4) &= -(-1)^k\nabla^{k+1}R^{\nabla}(3,...,k+2,2,k+3,1,k+4)\\
&- (-1)^k\, \nabla^{k+1}R^{\nabla}\,(3,...,k+2,1,2,k+3,k+4)\\
&+ (Q_k+V_k)(1,...,k+4)\,,
\end{aligned}$$

thanks to the alternating property of $R^{\nabla}$. Then the identity (2.3) follows from the differential Bianchi identity. □

## 2.2 Equivalent definitions of the tensor $S_k$

We start by noticing a few elementary equivalent definitions of the tensors $S_k$ introduced in the statement of corollary 1.3. (We remind that we assume more in general here that $\nabla$ is complex). For notations simplicity we will use the identification $\nabla^k S\equiv S^{(k)}$ for any tensor $S$. Using the identity (2.2) we infer the following equivalent definition for $S_k$.

**Definition 2.6.** *We define for all $k\geqslant -1$*

$$\Phi_k(1,...,k+4) \;:=\; -(-1)^k\, R^{(k+1)}(2,...,k+3,1,k+4)\,,$$

*We define also for all $k \geqslant 2$*

$$S_k \;:=\; \frac{i}{(k+1)!k!}\,\mathrm{Sym}_{2,\ldots,k+1}\,\theta_k\,,$$

$$\theta_k \;:=\; 2i^k\Phi_{k-3} + k!\sum_{p=2}^{k-2}(pS_p\wedge S_{k-p}) + \sum_{r=3}^{k-2}(r+1)!\sum_{p=2}^{r-1}\rho_{k,r,p}\,,$$

$$\rho_{k,r,p}(1,\ldots,k+1) \;:=\; (-i)^{k-1-r}\,(pS_p\wedge S_{r-p+1})^{(k-1-r)}(2,\ldots,k-r,1,k-r+1,\ldots,k+1)\,.$$

**Remark 2.7.** Notice that

$$\theta_4 \;=\; 2\Phi_1 + 4!2S_2\wedge S_2\,,$$

and

$$\mathrm{Sym}_{2,3,4,5}(S_2\wedge S_2) \;=\; \mathrm{Sym}_{2,3,4,5}\,\mu_4\,,$$

$$\mu_4(1,\ldots,5) \;:=\; -S_2(2,S_2(1,3,4),5)\,.$$

This implies

$$\mathrm{Sym}_{2,3,4,5}(S_2\wedge S_2) \;=\; \frac{1}{18}\,\mathrm{Sym}_{2,3,4,5}\,\tilde{\mu}_4\,,$$

$$\tilde{\mu}_4(1,\ldots,5) \;=\; R(2,R(1,3,4),5)\,.$$

**Remark 2.8.** Notice that

$$\theta_5 \;=\; 2i\Phi_2 + 5!2S_2\wedge S_3 + 5!3S_3\wedge S_2 + 4!\rho_{5,3,2}\,,$$

$$\rho_{5,3,2}(1,\ldots,6) \;=\; -2i(S_2\wedge S_2)'(2,1,3,4,5,6)\,,$$

and

$$\mathrm{Sym}_{2,3,4,5,6}\,\theta_5 \;=\; \mathrm{Sym}_{2,3,4,5,6}\,(2i\Phi_2 + \mu_5)\,,$$

$$\mu_5(1,\ldots,6) \;:=\; -5!2S_2(2,S_3(1,3,4,5),6) - 5!3S_3(2,S_2(1,3,4),5,6)$$

$$+\; 2i4!S_2'(2,3,S_2(1,4,5),6) + 2i4!S_2(3,S_2'(2,1,4,5),6)\,.$$

$$\mathrm{Sym}_{2,3,4,5,6}\,\mu_5 \;=\; \mathrm{Sym}_{2,3,4,5,6}\,\tilde{\mu}_5$$

$$\tilde{\mu}_5(1,\ldots,6) \;:=\; -6iR(2,R'(3,1,4,5),6) - 6iR'(2,3,R(1,4,5),6)\,.$$

**Remark 2.9.** We observe the equality

$$(S_p\wedge_1 S_{r-p+1})(1,k-r+1,\ldots,k+1)$$

$$:=\; S_p(1,S_{r-p+1}(k-r+1,\ldots,k-p+2),k-p+3,\ldots,k+1)$$

$$-\; S_p(k-r+1,S_{r-p+1}(1,k-r+2,\ldots,k-p+2),k-p+3,\ldots,k+1)\,.$$

Now, the fact that the tensors $\theta_k$ have at least one couple of alternating entries combined with the fact that the total symmetrization of such tensors vanishes implies the following equivalent definition of $S_k$.

**Definition 2.10.** ***Equivalent definition of $S_k$.*** *We define for all $k \geqslant -1$*

$$\Phi_k(1, ..., k+4) \;:=\; -(-1)^k R^{(k+1)}(2, ..., k+3, 1, k+4)\,,$$

*We define also for all $k \geqslant 2$*

$$S_k \;:=\; \frac{i}{(k+1)!k!} \operatorname{Sym}_{2,\dots,k+1} \theta_k\,,$$

$$\theta_k \;:=\; 2i^k \Phi_{k-3} + \mu_k$$

$$\mu_k \;:=\; -\sum_{r=3}^{k-1} (r+1)!(-i)^{k-1-r} \sum_{p=2}^{r-1} p \sum_{I \subseteq \{2,\dots,k-r\}} \rho^I_{k,r,p}\,,$$

*with*

$$\rho^I_{k,r,p}(1, ..., k+1)$$

$$:=\; S_p^{(|I|)}\big(I, k-r+1, S^{(|\complement I|)}_{r-p+1}(\complement I, 1, k-r+2, ..., k-p+2), k-p+3, ..., k+1\big)\,.$$

(Notice that in the case $r = k-1$ the set $I$ is empty). Using the elementary properties of the symmetrization operators we infer the following equivalent definition.

**Definition 2.11.** ***Equivalent definition of $S_k$.*** *We define for all $k \geqslant -1$*

$$\Phi_k(1, ..., k+4) \;:=\; -(-1)^k R^{(k+1)}(2, ..., k+3, 1, k+4)\,,$$

*We define also for all $k \geqslant 2$*

$$\boldsymbol{\theta}_k \;:=\; 2i^k \Phi_{k-3} + \boldsymbol{\mu}_k\,,$$

$$\boldsymbol{\mu}_k \;:=\; \sum_{r=3}^{k-1} \sum_{p=2}^{r-1} \frac{(r+1)!(-i)^{k-1-r}}{(r-p+2)!\,(p+1)!} \sum_{I \subseteq \{2,\dots,k-r\}} \sum_{j=k-p+2}^{k+1} \rho^I_{k,r,p,j}\,,$$

*with*

$$\rho^I_{k,r,p,j}(1, ..., k+1)$$

$$:=\; \boldsymbol{\theta}_p^{(|I|)}\big(I, k-r+1, k-p+3, ..., j, \boldsymbol{\theta}^{(|\complement I|)}_{r-p+1}(\complement I, 1, k-r+2, ..., k-p+2), j+1, ..., k+1\big)\,.$$

*We define also for all $k \geqslant 2$*

$$S_k \;:=\; \frac{i}{(k+1)!k!} \operatorname{Sym}_{2,\dots,k+1} \boldsymbol{\theta}_k\,.$$

**Proposition 2.12.** *For any integer $k \geqslant 2$ hold the identity*

$$S_k \;=\; \frac{i}{(k+1)!k!} \operatorname{Sym}_{2,\dots,k+1} \tilde{\boldsymbol{\theta}}_k\,,$$

*with*

$$\tilde{\boldsymbol{\theta}}_k \;=\; \sum_{D \vDash_0 k-2l_D} C_D\, R^D\,,$$

*where the coefficients $C_D$ are given by $C_D := 2\,(-i)^k$ for $l_D = 1$ and by the recursive formula*

$$C_D = -\sum_{\substack{1\leqslant h\leqslant l_D-1\\ 0\leqslant A\leqslant D\\ l_A=l_D}} (-i)^{|D-A|}\binom{|D-A|}{D-A}\frac{\|A\|!\,C_{A_h^-}\,C_{A_h^+}}{(\|A_h^-\|+1)!(\|A_h^+\|+1)!}\,, \tag{2.4}$$

*for $l_D \geqslant 2$, where $A_h^- := (a_1, ..., a_h)$ and $A_h^+ := (a_{h+1}, ..., a_{l_D})$ for any $A = (a_1, ..., a_{l_D})$.*

**Proof.** We assume by induction on $k \geqslant 3$ that

$$S_p \;=\; \frac{i}{(p+1)!\,p!}\,\mathrm{Sym}_{2,...,p+1}\,\tilde{\boldsymbol{\theta}}_p\,,$$

for any $p = 2, ...., k-1$, with $\tilde{\boldsymbol{\theta}}_p$ under the form

$$\tilde{\boldsymbol{\theta}}_p \;=\; \sum_{A\vDash_0 p-2l_A} C_A\,R^A\,,$$

for some coefficients $C_A$, with obviously $\tilde{\boldsymbol{\theta}}_2 = 2R$ and $C_A = 2\,(-i)^p$ in the case $l_A = 1$, thanks to the equivalent definition (2.11) of $S_p$. We will show that $\tilde{\boldsymbol{\theta}}_k$ writes under the form claimed in the statement of the proposition 2.12 with the coefficients $C_D$ given by the recursive formula (2.4). For this purpose we notice first that the straightforward argument showing the equivalent definition (2.11) of $S_k$ implies for all $k \geqslant 3$ the identity

$$S_k \;=\; \frac{i}{(k+1)!\,k!}\,\mathrm{Sym}_{2,...,k+1}\,\hat{\boldsymbol{\theta}}_k\,,$$

with

$$\hat{\boldsymbol{\theta}}_k \;:=\; 2\,i^k\Phi_{k-3} + \hat{\boldsymbol{\mu}}_k\,,$$

$$\hat{\boldsymbol{\mu}}_k \;:=\; \sum_{r=3}^{k-1}\sum_{p=2}^{r-1}\frac{(r+1)!(-i)^{k-1-r}}{(r-p+2)!\,(p+1)!}\sum_{I\subseteq\{2,...,k-r\}}\;\sum_{j=k-p+2}^{k+1}\hat{\rho}^I_{k,r,p,j}\,,$$

and with

$$\hat{\rho}^I_{k,r,p,j}(1,...,k+1)$$

$$:=\; \tilde{\boldsymbol{\theta}}_p^{(|I|)}\big(I, k-r+1, k-p+3, ..., j,\, \tilde{\boldsymbol{\theta}}_{r-p+1}^{(|\complement I|)}(\complement I, 1, k-r+2, ..., k-p+2), j+1, ..., k+1\big)\,.$$

For any integer $a \geqslant 1$, we denote by $[a] := \{1, ..., a\}$ and we set $[0] := \emptyset$. For any integer $b \geqslant 1$, we denote by $\mathrm{Map}(a,b)$ the set of maps $f\colon [a] \longrightarrow [b]$.

We observe now that for any $A \vDash_0 p - 2l_A$ and any integer $q \geqslant 0$, the Leibniz identity implies

$$(R^A)^{(q)}(1, ..., q+p+1) \;=\; \sum_{f\in\mathrm{Map}(q,l_A)}\left[\prod_{j=1}^{l_A} R^{A,q}_{f,j}\right]\cdot(q+1)\,,$$

with

$$R^{A,q}_{f,j} \;:=\; R^{(a_j+|f^{-1}(j)|)}\big(f^{-1}(j)\,, \|A\|_{j-1}+q+2, ..., \|A\|_j+q+1, \bullet, p-j+q+2\big)\,.$$

We write now

$$\tilde{\boldsymbol{\theta}}_p^{(|I|)} \;=\; \sum_{A \vDash_0 p-2l_A} C_A\,(R^A)^{(|I|)}\,,$$

and we notice that the shape of $(R^A)^{(|I|)}$ shows that $\hat{\rho}^I_{k,r,p,k-h+1}$, with $h := l_A$, are the only terms with non vanishing symmetrization of the variables 2, ..., $k+1$. Indeed we consider the factor

$$R^{(a_h+|f^{-1}(h)|)}(f^{-1}(h), p-h+1-a_h\,, ..., p-h+1, 1, p-h+2)\,, \tag{2.5}$$

with $f \in \mathrm{Map}(I, h)$ in the expression of $(R^A)^{(|I|)}(I, 1, ..., p+1)$ and we perform the change of variables

$$\begin{array}{cccccccc}
1 & 2 & \ldots & j-k+p-1 & j-k+p & j-k+p & \ldots & p+1\\
\downarrow & \downarrow & & \downarrow & \downarrow & \downarrow & & \downarrow\\
k-r+1 & k-p+3 & \ldots & j & \tilde{\boldsymbol{\theta}}^{(|\complement I|)}_{r-p+1}(\complement I, 1, k-r+2, ..., k-p+2) & j+1 & \ldots & k+1
\end{array}$$

in the factor (2.5) with $j \in \{k-p+2, ..., k+1\}$. (We remind that we ignore the standard increasing notation $a, ..., b$ when $a > b$.) The only case when the symmetrization of the variables $2, ..., k+1$ does not annihilate the factor (2.5), is when the index $j \in \{k-p+2, ..., k+1\}$ satisfies the equality

$$j-k+p \;=\; p-h+1,$$

i.e. only when $j = k-h+1$. This shows the required statement about $\hat{\rho}^I_{k,r,p,k-h+1}$.

We infer that for all $k \geqslant 3$ hold the equality

$$S_k \;=\; \frac{i}{(k+1)!k!}\,\mathrm{Sym}_{2,...,k+1}\,\check{\boldsymbol{\theta}}_k\,,$$

with

$$\check{\boldsymbol{\theta}}_k \;:=\; 2i^k\Phi_{k-3} + \check{\boldsymbol{\mu}}_k\,,$$

with

$$\check{\boldsymbol{\mu}}_k \;:=\; \sum_{\substack{3 \leqslant r \leqslant k-1\\ 2 \leqslant p \leqslant r-1}} \frac{(r+1)!(-i)^{k-1-r}}{(r-p+2)!\,(p+1)!} \sum_{\substack{A \vDash_0 p-2l_A\\ B \vDash_0 r-p+1-2l_B\\ I \subseteq \{2,...,k-r\}}} C_A\,C_B\,\check{\rho}^{I,A,B}_{k,r,p}\,,$$

with

$$\begin{aligned}
&\check{\rho}^{I,A,B}_{k,r,p}(1, ..., k+1)\\
:=\;& (R^A)^{(|I|)}(I, k_r+1, k_p+3, ..., k_{l_A}+1, (R^B)^{(|\complement I|)}(\complement I, 1, k_r+2, ..., k_p+2), k_{l_A}+2, ..., k+1)\,,
\end{aligned}$$

and with $k_s := k-s$. Using the expression

$$R^A(x, y, ..., p+y-1) \;=\; \left[\prod_{j=1}^{l_A} R^{(a_j)}\big(\|A\|_{j-1}+y\,, ..., \|A\|_j+y-1, \bullet, p-j+y\big)\right]\cdot x\,,$$

for any integers $x, y > 0$ and using the change of variables

$$\begin{array}{ccccccc} k_p+3 & \dots & k_{l_A}+1 & k_{l_A}+2 & k_{l_A}+3 & \dots & k+2 \\ \downarrow & & \downarrow & \downarrow & \downarrow & & \downarrow \\ k_p+3 & \dots & k_{l_A}+1 & R^B(1, k_r+2, \dots, k_p+2) & k_{l_A}+2 & \dots & k+1 \end{array}$$

we infer the equality

$$R^A(k_r+1, k_p+3, \dots, k_{l_A}+1, R^B(1, k_r+2, \dots, k_p+2), k_{l_A}+2, \dots, k+1)$$

$$= \left[\prod_{j=1}^{l_A-1} R^{(a_j)}\big(\|A\|_{j-1}+k_p+3, \dots, \|A\|_j+k_p+2, \bullet, k_j+2\big)\right] \cdot$$

$$\cdot \quad R^{(a_{l_A})}\big(k_{l_A}+2-a_{l_A}, \dots, k_{l_A}+1, R^B(1, k_r+2, \dots, k_p+2), k_r+1, k_{l_A}+2\big).$$

Indeed notice that $|A| + l_A = p - l_A = k_{l_A} + 2 - (k_p + 2)$ and $p - j + (k_p + 3) - 1 = k_j + 2$. Expanding $R^B(1, k_r+2, \dots, k_p+2)$ we obtain

$$R^A(k_r+1, k_p+3, \dots, k_{l_A}+1, R^B(1, k_r+2, \dots, k_p+2), k_{l_A}+2, \dots, k+1)$$

$$= \left[\prod_{j=1}^{l_A-1} R^{(a_j)}\big(\|A\|_{j-1}+k_p+3, \dots, \|A\|_j+k_p+2, \bullet, k_j+2\big)\right] \cdot$$

$$\cdot \quad R^{(a_{l_A})}(k_{l_A}+2-a_{l_A}, \dots, k_{l_A}+1, \bullet, k_r+1, k_{l_A}+2) \cdot$$

$$\cdot \quad \left[\prod_{t=1}^{l_B} R^{(b_t)}\big(\|B\|_{t-1}+k_r+2, \dots, \|B\|_t+k_r+1, \bullet, k_p-t+3\big)\right] \cdot 1,$$

since $\|B\| - t + k_r + 2 = k_p - t + 3$.

For any integer $a \geqslant 1$, we denote by $[a]_1 := \{2, \dots, a+1\}$ and we set $[0]_1 := \emptyset$. For any integer $b \geqslant 1$, we denote by $\mathrm{Map}([a]_1, b)$ the set of maps $f \colon [a]_1 \longrightarrow [b]$. Using the Leibniz formula and the alternating property of $R$ we infer the identity

$$\check{\boldsymbol{\mu}}_k = \sum_{\substack{3 \leqslant r \leqslant k-1 \\ 2 \leqslant p \leqslant r-1}} \frac{(r+1)!(-i)^{k-1-r}}{(r-p+2)!\,(p+1)!} \sum_{\substack{A \vDash_0 p-2l_A \\ B \vDash_0 r-p+1-2l_B \\ f \in \mathrm{Map}([k_r-1]_1, l_A+l_B)}} C_A\, C_B\, \check{\rho}_{k,r,p}^{f,A,B}, \tag{2.6}$$

with

$$- \quad \check{\rho}_{k,r,p}^{f,A,B}(1, \dots, k+1)$$

$$:= \left[\prod_{j=1}^{l_A-1} R^{(a_j+|f^{-1}(j)|)}(f^{-1}(j), \|A\|_{j-1}+k_p+3, \dots, \|A\|_j+k_p+2, \bullet, k_j+2)\right] \cdot$$

$$\cdot \quad R^{(a_{l_A}+|f^{-1}(l_A)|)}\big(f^{-1}(l_A), k_{l_A}+2-a_{l_A}, \dots, k_{l_A}+1, k_r+1, \bullet, k_{l_A}+2\big) \cdot$$

$$\cdot \quad \left[\prod_{t=1}^{l_B} R^{(b_t+|f^{-1}(l_A+t)|)}\big(f^{-1}(l_A+t), \|B\|_{t-1}+k_r+2, \dots, \|B\|_t+k_r+1, \bullet, k_p-t+3\big)\right] \cdot 1.$$

We observe the equality

$$-\operatorname{Sym}_{2,\ldots,k+1}\check{\rho}_{k,r,p}^{f,A,B}=\operatorname{Sym}_{2,\ldots,k+1}R^D\,, \tag{2.7}$$

with $l_D:=l_A+l_B$ and with $D\vDash_0 k-2l_D$ given by

$$d_j:=\begin{cases} a_j+|f^{-1}(j)|\,, & j=1,\ldots,l_A\,,\\ \\ b_{j-l_A}+|f^{-1}(j)|\,, & j=l_A+1,\ldots,l_D\,.\end{cases} \tag{2.8}$$

We deduce that $\tilde{\boldsymbol{\theta}}_k$ writes under the form claimed in the statement of the proposition 2.12.

We denote by abuse of notations $D-A-B\geqslant 0$ when $d_j-a_j\geqslant 0$ for all $j=0,\ldots,l_A$ and when $d_j-b_{j-l_A}\geqslant 0$ for all $j=l_A+1,\ldots,l_D$.

Given $A$, $B$, $D$ as before we denote by $\operatorname{Map}_{A,B}^D([k-r-1]_1,l_D)$ the sub-set of $\operatorname{Map}([k-r-1]_1,l_D)$ given by the elements $f$ which satisfy condition (2.8). We denote also by

$$M_{A,B,k,r}^D\;:=\;|\operatorname{Map}_{A,B}^D([k-r-1]_1,l_D)|\,.$$

Using (2.6) and (2.7), we infer the recursive formula

$$C_D\;=\;\sum_{\substack{3\leqslant r\leqslant k-1\\2\leqslant p\leqslant r-1}}\frac{(r+1)!(-i)^{k+1-r}}{(r-p+2)!\,(p+1)!}\sum_{\substack{A\vDash_0 p-2l_A\\B\vDash_0 r-p+1-2l_B\\l_A+l_B=l_D\\D-A-B\geqslant 0}}C_A\,C_B\,M_{A,B,k,r}^D\,,$$

We remind that in the case $l_D=1$ hold the formula $C_D=2\,(-i)^k$ for all $k\geqslant 2$.

**Remark 2.13.** Let $A\vDash_q a$ and consider the sub-set $\operatorname{Map}_A(a,q)$ of $\operatorname{Map}(a,q)$ given by the elements $f$ such that $|f^{-1}(j)|=a_j$, for all $j=1,\ldots,q$. Then

$$|\operatorname{Map}_A(a,q)|\;=\;\binom{a}{A}.$$

The fact that in our set-up $|D|-|A|-|B|=k-r-1$, allows to apply the previous remark to the set $\operatorname{Map}_{A,B}^D([k-r-1]_1,l_D)$. We infer the formula

$$M_{A,B,k,r}^D\;=\;(k-r-1)!\prod_{j=1}^{l_A}\frac{1}{(d_j-a_j)!}\prod_{j=l_A+1}^{l_D}\frac{1}{(d_j-b_{j-l_A})!}\,.$$

We conclude the explicit recursive formula

$$C_D\;=\;\sum_{\substack{3\leqslant r\leqslant k-1\\2\leqslant p\leqslant r-1\\A\vDash_0 p-2l_A\\B\vDash_0 r-p+1-2l_B\\l_A+l_B=l_D\\D-A-B\geqslant 0}}\frac{(-i)^{k+1-r}\,(k-r-1)!(r+1)!C_A\,C_B}{(r-p+2)!\,(p+1)!\prod_{j=1}^{l_A}(d_j-a_j)!\prod_{j=l_A+1}^{l_D}(d_j-b_{j-l_A})!}\,,$$

which rewrites as (2.4). This concludes the proof of proposition 2.12. □

**Example 2.14.** In this example we set for simplicity $R'$: $=R^{(1)}$ and $R'' := R^{(2)}$. For $k=6$ we have the expression

$$\begin{aligned}
\tilde{\boldsymbol{\theta}}_6(1,...,7) &= -2\Phi_3(1,...,7)\\
&+ C_{2,0}\,R''(2,3,4,R(5,1,6),7)+C_{0,2}\,R(2,R''(3,4,5,1,6),7)\\
&+ C_{1,1}\,R'(2,3,R'(4,5,1,6),7)+C_{0,0,0}\,R(2,R(3,R(4,1,5),6),7)\,,\\
C_{2,0} &= 10,\\
C_{0,2} &= 10,\\
C_{1,1} &= 17,\\
C_{0,0,0} &= -\frac{32}{3}\,.
\end{aligned}$$

# 3 Expliciting the recursive formula for $C_D$

**Lemma 3.1.** *For any integer $k\geqslant 3$ and any $D\vDash_0 k-2l_D$ with $l_D\geqslant 1$, hold the recursive formula*

$$C_D \;=\; -i\sum_{p\in\operatorname{Supp}D} C_{D-1_p} - k!\sum_{h=1}^{l_D-1}\frac{C_{D_h^-}\,C_{D_h^+}}{(\|D_h^-\|+1)!(\|D_h^+\|+1)!}\,.$$

**Proof.** In the case $l_D=1$ we have $C_{k-2}=-iC_{k-3}$, for $k\geqslant 3$ by definition. This coincides with what provides the formula in the lemma under consideration. We consider now the case $l_D\geqslant 2$ and we use the identifications $D\equiv(d_1,...,d_{l_D})$, $A\equiv(a_1,...,a_{l_D})$. Using the recursive formula (2.4) we write

$$\begin{aligned}
&-i\sum_{p\in\operatorname{Supp}D} C_{D-1_p}\\
&= -\sum_{\substack{p\in\operatorname{Supp}D\\ 0\leqslant A\leqslant D-1_p\\ l_A=l_D}} (-i)^{|D|-|A|}\,\frac{d_p-a_p}{|D|-|A|}\binom{|D|-|A|}{D-A}F_A\,,
\end{aligned}$$

with

$$F_A \;:=\; \|A\|!\sum_{h=1}^{l_D-1}\frac{C_{A_h^-}\,C_{A_h^+}}{(\|A_h^-\|+1)!(\|A_h^+\|+1)!}\,.$$

Let $\mathbb{I}\equiv(1,\,...,\,1)\in\mathbb{Z}^{l_D}$ and let $S_D:=|\operatorname{Supp}\,D|$. We denote by $\mathbb{I}_D\in\mathbb{Z}^{l_D}$ the vector such that $\operatorname{Supp}\mathbb{I}_D=\operatorname{Supp}D$ and $0\leqslant\mathbb{I}_D\leqslant\mathbb{I}$. Then

$$\begin{aligned}
&-i\sum_{p\in\operatorname{Supp}D} C_{D-1_p}\\
&= -\sum_{\substack{0\leqslant A\leqslant D-\mathbb{I}_D\\ l_A=l_D}} (-i)^{|D|-|A|}\left[\sum_{\substack{1\leqslant p\leqslant l_D\\ p\in\operatorname{Supp}D}}\frac{d_p-a_p}{|D|-|A|}\right]\binom{|D|-|A|}{D-A}F_A
\end{aligned}$$

$$- \sum_{\substack{p \in \mathrm{Supp} D \\ D - \mathbb{I}_D \leqslant A \leqslant D - 1_p \\ l_A = l_D \\ |D| - S_D < |A|}} \frac{(-i)^{|D|-|A|}}{|D|-|A|} \binom{|D|-|A|}{D-A} F_A .$$

We notice indeed that the conditions $D - \mathbb{I}_D \leqslant A \leqslant D - 1_p$ , $|D| - S_D < |A|$, when they are not empty, imply $a_p = d_p - 1$. We infer

$$\begin{aligned}
& -i \sum_{p \in \mathrm{Supp} D} C_{D-1_p} \\
= & - \sum_{\substack{0 \leqslant A \leqslant D - \mathbb{I}_D \\ l_A = l_D}} (-i)^{|D|-|A|} \binom{|D|-|A|}{D-A} F_A \\
- & \sum_{\substack{D - \mathbb{I}_D \leqslant A \leqslant D \\ l_A = l_D \\ |D| - S_D < |A| < |D|}} |\{p : a_p = d_p - 1\}| \frac{(-i)^{|D|-|A|}}{|D|-|A|} \binom{|D|-|A|}{D-A} F_A \\
= & - \sum_{\substack{0 \leqslant A \leqslant D \\ l_A = l_D \\ |A| < |D|}} (-i)^{|D|-|A|} \binom{|D|-|A|}{D-A} F_A ,
\end{aligned}$$

since $|\{p : a_p = d_p - 1\}| = |D| - |A|$. Then the conclusion follows from the formula (2.4). □

**Definition 3.2.** *Given $t \in \mathbb{R}^{\mathbb{Z}_{>0}}$, $s \in \mathbb{R}$ and $A \in \mathbb{Z}^p_{\geqslant 0}$ , $B \in \mathbb{Z}^q_{\geqslant 0}$, with $p, q \in \mathbb{Z}_{>0}$ we define the concatenation product*

$$(t^A s^p) * (t^B s^q) \ := \ t_1^{a_1} \cdots t_p^{a_p} t_{p+1}^{b_1} \cdots t_{p+q}^{b_q} s^{p+q} ,$$

*and we extend it by linearity.*

**Definition 3.3.** *Given $t \in \mathbb{R}^{\mathbb{Z}_{>0}}$, $s \in \mathbb{R}$, $A \in \mathbb{Z}^p_{\geqslant 0}$ and $j \in \mathbb{Z}_{>0}$ we define the contraction product*

$$t_j \neg (t^A s^p) \ := \ t^{A + 1_j} s^p ,$$

*where $A + 1_j := A + 1_{p,j}$, if $j \leqslant p$, with $1_{p,j}$ the vector of length $p$ with vanishing components except the $j$-th one which has value 1 and where $A + 1_j := A$, if $j > p$. We extend the contraction product by linearity.*

**Lemma 3.4.** *Let $t \in \mathbb{R}^{\mathbb{Z}_{>0}}$, $s \in \mathbb{R}$, and $x \in \mathbb{R}$. Then the function*

$$u(x, t, s) \ := \ \sum_{D \geqslant 0} C_D t^D s^{l_D} \frac{x^{\|D\|+1}}{(\|D\|+1)!} ,$$

*is the unique solution of the Riccati-type ODE*

$$\partial_x u + \frac{u * u}{x^2} + i \left( \sum_{j \geqslant 1} t_j \right) \neg u + s x^2 \ = \ 0,$$

*with the initial condition* $u(0,\cdot,\cdot)\equiv 0$.

**Proof.** We write first

$$\begin{aligned}
\partial_x u(x,t,s) &= \sum_{D\geqslant 0} C_D\, t^D\, s^{l_D}\, \frac{x^{\|D\|}}{\|D\|!} \\
&= \sum_{\substack{D\geqslant 0\\ |D|+2l_D\geqslant 3}} C_D\, t^D s^{l_D}\, \frac{x^{\|D\|}}{\|D\|!} + C_{2,\mathbf{0}}\, s\, \frac{x^2}{2} \\
&= \sum_{\substack{D\geqslant 0\\ |D|+2l_D\geqslant 3}} C_D\, t^D s^{l_D}\, \frac{x^{\|D\|}}{\|D\|!} - s x^2 .
\end{aligned}$$

On the other hand

$$\begin{aligned}
(u*u)(x,t,s) &= \left[\sum_{A\geqslant 0} C_A\, t^A\, s^{l_A}\, \frac{x^{\|A\|+1}}{(\|A\|+1)!}\right] * \left[\sum_{B\geqslant 0} C_B\, t^B\, s^{l_B}\, \frac{x^{\|B\|+1}}{(\|B\|+1)!}\right] \\
&= x^2 \sum_{A,B\geqslant 0} C_A\, C_B\, (t^A\, s^{l_A}) * (t^B\, s^{l_B})\, \frac{x^{\|A\|+\|B\|}}{(\|A\|+1)!(\|B\|+1)!} \\
&= x^2 \sum_{D\geqslant 0} t^D\, s^{l_D} \sum_{h=1}^{l_D-1} \frac{C_{D_h^-}\, C_{D_h^+}\, x^{\|D\|}}{(\|D_h^-\|+1)!(\|D_h^+\|+1)!} \\
&= x^2 \sum_{\substack{D\geqslant 0\\ |D|+2l_D\geqslant 3}} t^D\, s^{l_D} \sum_{h=1}^{l_D-1} \frac{C_{D_h^-}\, C_{D_h^+}\, x^{\|D\|}}{(\|D_h^-\|+1)!(\|D_h^+\|+1)!},
\end{aligned}$$

and

$$\begin{aligned}
\left[\left(\sum_{j\geqslant 1} t_j\right)\neg u\right](x,t,s) &= \sum_{j\geqslant 1} (t_j \neg u)(x,t,s) \\
&= \sum_{j\geqslant 1} \sum_{D\geqslant 0} C_D\, t^{D+1_j} s^{l_D}\, \frac{x^{\|D\|+1}}{(\|D\|+1)!} \\
&= \sum_{D\geqslant 0} \sum_{j=1}^{l_D} C_D\, t^{D+1_j}\, s^{l_D}\, \frac{x^{\|D\|+1}}{(\|D\|+1)!} \\
&= \sum_{D\geqslant 0} t^D\, s^{l_D}\, \frac{x^{\|D\|}}{\|D\|!} \left(\sum_{p\in \mathrm{Supp}D} C_{D-1_p}\right) \\
&= \sum_{\substack{D\geqslant 0\\ |D|+2l_D\geqslant 3}} t^D\, s^{l_D}\, \frac{x^{\|D\|}}{\|D\|!} \left(\sum_{p\in \mathrm{Supp}D} C_{D-1_p}\right).
\end{aligned}$$

Then lemma (3.1) implies that the function $u$ satisfies the singular Riccati's EDO. The uniqueness statement follows from the following corollary. □

**Corollary 3.5.** *Let $t\in\mathbb{R}^{\mathbb{Z}_{>0}}$, $s\in\mathbb{R}$, and $x\in\mathbb{R}$ then the function*

$$u(x,t,s) \;=\; \sum_{l\geqslant 1} u_l(x,t_1,\ldots,t_l)\, s^l\,,$$

$$u_l(x,t_1,\ldots,t_l) \;:=\; \sum_{\substack{D\geqslant 0\\ l_D=l}} C_D\, t^D \frac{x^{\|D\|+1}}{(\|D\|+1)!}\,,$$

*satisfies the recursive system of linear ODE's*

$$\begin{cases} \partial_x u_1(x,t_1)+i t_1 u_1(x,t_1)+x^2=0\,, \\ \partial_x u_l(x,t_1,\ldots,t_l)+i u_l(x,t_1,\ldots,t_l)\displaystyle\sum_{j=1}^{l} t_j+\frac{1}{x^2}\sum_{p=1}^{l-1} u_p(x,t_1,\ldots,t_p)\, u_{l-p}(x,t_{p+1},\ldots,t_l) \;=\; 0\,, \\ \forall l\geqslant 2\,. \end{cases} \tag{3.1}$$

*with initial conditions $u_l(0,t_1,\ldots,t_l)=0$, for all $l\geqslant 1$ and all $(t_1,\ldots,t_l)\in\mathbb{R}^l$.*

**Proof.** According to the proof of lemma 3.4 we can write

$$\begin{aligned}
&(u*u)(x,t,s)\\
&= x^2\sum_{A,B\geqslant 0} C_A\, C_B\, (t^A s^{l_A})*(t^B s^{l_B})\,\frac{x^{\|A\|+\|B\|}}{(\|A\|+1)!(\|B\|+1)!}\\
&= x^2\sum_{l\geqslant 2} s^l\sum_{p=1}^{l-1}\sum_{\substack{A\geqslant 0\\ l_A=p}}\sum_{\substack{B\geqslant 0\\ l_B=l-p}} \frac{C_A\, C_B\, t_1^{a_1}\cdots t_p^{a_p}\, t_{p+1}^{b_1}\cdots t_l^{b_{l-p}}\, x^{\|A\|+\|B\|}}{(\|A\|+1)!(\|B\|+1)!}\\
&= \sum_{l\geqslant 2} s^l\sum_{p=1}^{l-1}\left[\sum_{\substack{A\geqslant 0\\ l_A=p}} \frac{C_A\, t_1^{a_1}\cdots t_p^{a_p}\, x^{\|A\|+1}}{(\|A\|+1)!}\right]\left[\sum_{\substack{B\geqslant 0\\ l_B=l-p}} \frac{C_B\, t_{p+1}^{b_1}\cdots t_l^{b_{l-p}}\, x^{\|B\|+1}}{(\|B\|+1)!}\right]\\
&= \sum_{l\geqslant 2} s^l\sum_{p=1}^{l-1} u_p(x,t_1,\ldots,t_p)\, u_{l-p}(x,t_{p+1},\ldots,t_l)\,.
\end{aligned}$$

On the other hand still according to the proof of lemma 3.4 we can write

$$\begin{aligned}
\left[\left(\sum_{j\geqslant 1} t_j\right)\neg u\right](x,t,s) &= \sum_{D\geqslant 0}\sum_{j=1}^{l_D} C_D\, t^{D+1_j}\, s^{l_D}\frac{x^{\|D\|+1}}{(\|D\|+1)!}\\
&= \sum_{l\geqslant 1} s^l\sum_{\substack{D\geqslant 0\\ l_D=l}}\sum_{j=1}^{l} C_D\, t^D t_j\frac{x^{\|D\|+1}}{(\|D\|+1)!}\\
&= \sum_{l\geqslant 1} s^l\left(\sum_{j=1}^{l} t_j\right) u_l(x,t_1,\ldots,t_l)\,.
\end{aligned}$$

The required conclusion follows from the Riccati-type equation in the statement of lemma 3.4. □

More in general we introduce the concatenation product for functions of the type

$$u(x,t,s) = \sum_{l\geqslant 1} u_l(x,t_1,...,t_l)\, s^l ,$$

$$v(x,t,s) = \sum_{l\geqslant 1} v_l(x,t_1,...,t_l)\, s^l ,$$

as

$$(u*v)(x,t,s) = \sum_{l\geqslant 2} s^l \sum_{p=1}^{l-1} u_p(x,t_1,...,t_p)\, v_{l-p}(x,t_{p+1},...,t_l)\,.$$

The concatenation product is associative.

**Lemma 3.6.** *The solution $u$ of the system (3.1) is equivalent to the solution*

$$v(x,t,s) = \sum_{l\geqslant 1} v_l(x,t_1,...,t_l)\, s^l ,$$

*of the Riccati-type ODE*

$$\partial_x v - \frac{i}{x^2}\, v*v + i x^2\, e^{ixt_1} s = 0, \tag{3.2}$$

*via the identification*

$$u_l(x,t_1,...,t_l) = -i\, e^{-ix\sum_{j=1}^l t_j}\, v_l(x,t_1,...,t_l)\,.$$

**Proof.** We write

$$\partial_x u_l(x,t_1,...,t_l) = -i\, e^{-ix\sum_{j=1}^l t_j}\, \partial_x v_l(x,t_1,...,t_l) - e^{-ix\sum_{j=1}^l t_j}\, v_l(x,t_1,...,t_l) \sum_{j=1}^{l} t_j ,$$

$$i\, u_l(x,t_1,...,t_l) \sum_{j=1}^{l} t_j = e^{-ix\sum_{j=1}^l t_j}\, v_l(x,t_1,...,t_l) \sum_{j=1}^{l} t_j ,$$

and

$$\sum_{p=1}^{l-1} u_p(x,t_1,...,t_p)\, u_{l-p}(x,t_{p+1},...,t_l)$$

$$= -\, e^{-ix\sum_{j=1}^l t_j} \sum_{p=1}^{l-1} v_p(x,t_1,...,t_p)\, v_{l-p}(x,t_{p+1},...,t_l)\,.$$

We infer that the system (3.1) rewrites as

$$\begin{cases} \partial_x v_1(x,t_1) + i x^2\, e^{ixt_1} = 0, \\ \partial_x v_l(x,t_1,...,t_l) - \dfrac{i}{x^2} \displaystyle\sum_{p=1}^{l-1} v_p(x,t_1,...,t_p)\, v_{l-p}(x,t_{p+1},...,t_l) = 0, \\ \forall l \geqslant 2, \end{cases} \tag{3.3}$$

which is equivalent to the equation (3.2). □

## 3.1 General facts about Riccati-type equations

In this subsection we denote for notations simplicity $\dot{y}(x,t,s) \equiv \partial_x y(x,t,s)$.

**Lemma 3.7.** *A solution*

$$y(x,t,s) \;=\; \sum_{l\geqslant 1} y(x,t_1,...,t_l)s^l\,,$$

*of the Ricati-type equation*

$$\dot{y} \;=\; a_2(x)y*y+a_1(x)y+a_0(x,t,s),$$

*with $a_2$ non zero and once differentiable and with*

$$a_0(x,t,s) \;=\; \sum_{l\geqslant 1} a_{0,l}(x,t_1,...,t_l)s^l\,,$$

*is equivalent to a solution*

$$u(x,t,s) \;=\; \sum_{l\geqslant 1} u(x,t_1,...,t_l)s^l\,,$$

*of the linear second order ODE*

$$\ddot{u}-R(x)\,\dot{u}+H*u \;=\; 0,$$

$$R \;:=\; a_1+\frac{\dot{a}_2}{a_2}\,,$$

$$H \;:=\; a_2\,a_0\,,$$

*via the identification $\dot{u}+a_2\,y*u=0$, with $u_1(x,t_1)\equiv f(t_1)$, for some non vanishing function and*

$$y_l(x,t_1,...,t_l)=\sum_{\lambda\vDash l}\,(-1)^{l_\lambda}\frac{\dot{u}_{\lambda_1+1}(x,t_1,...,t_{\lambda_1+1})}{a_2(x)f(t_{\lambda_1+1})}\prod_{j=2}^{l_\lambda}\,\frac{u_{\lambda_j+1}(x,t_{|\lambda|_{j-1}+1},...,t_{|\lambda|_j+1})}{f(t_{|\lambda|_j+1})}\,, \tag{3.4}$$

*for all $l\geqslant 1$, with $|\lambda|_j:=\sum_{p=1}^{j}\lambda_p$.*

**Proof.** The function $v:=a_2\,y$ satisfies Riccati-type equation

$$\dot{v} \;=\; v*v+R(x)\,v+H(x,t,s)\,.$$

Indeed

$$\begin{aligned}\dot{v} \;&=\; \dot{a}_2\,y+a_2\,\dot{y}\\ &=\; \frac{\dot{a}_2}{a_2}v+a_2\big[a_0+a_1\,y+a_2\,y*y\big]\\ &=\; a_0\,a_2+\left(a_1+\frac{\dot{a}_2}{a_2}\right)v+v*v.\end{aligned}$$

We show now that if $\dot{u}+v*u=0$ then the function $u$ satisfies the second order linear ODE in the statement. Indeed

$$\begin{aligned}\ddot{u} \;&=\; -\dot{v}*u-v*\dot{u}\\ &=\; -\dot{v}*u+v*(v*u)\end{aligned}$$

$$
\begin{aligned}
&= (-\dot{v} + v * v) * u \\
&= -[R(x)\, v + H] * u \\
&= R(x)\, \dot{u} - H * u.
\end{aligned}
$$

We notice now that the identification $\dot{u} + v*u = 0$ implies $\dot{u}_1 \equiv 0$, i.e. $u_1(x, t_1) \equiv f(t_1)$ and

$$y_l(x, t_1, ..., t_l) = -\frac{\dot{u}_{l+1}(x, t_1, ..., t_{l+1})}{a_2(x) f(t_{l+1})} - \sum_{p=1}^{l-1} y_p(x, t_1, ..., t_p) \frac{u_{l+1-p}(x, t_{p+1}, ..., t_{l+1})}{f(t_{l+1})}, \tag{3.5}$$

for all $l \geqslant 1$, which implies the explicit formula (3.4). Indeed if we set

$$
\begin{aligned}
A(x, t, s) &:= \sum_{l \geqslant 1} A_l(x, t_1, ..., t_l)\, s^l, \\
B(x, t, s) &:= \sum_{l \geqslant 1} B_l(x, t_1, ..., t_l)\, s^l, \\
A_l(x, t_1, ..., t_l) &:= \frac{u_{l+1}(x, t_1, ..., t_{l+1})}{f(t_{l+1})}, \\
B_l(x, t_1, ..., t_l) &:= \frac{\dot{u}_{l+1}(x, t_1, ..., t_{l+1})}{a_2(x) f(t_{l+1})},
\end{aligned}
$$

(notice that the definitions are well posed) then the relation (3.5) writes as $y = -B - y * A$. Thus

$$y = -B * \sum_{k \geqslant 0} (-1)^k A^{*,k},$$

where we denote by $A^{*,k}$ the $k$-th power with respect to the concatenation product $*$. (By convention $A^{*,0} := 1$.) In other terms

$$y_l(x, t_1, ..., t_l) = \sum_{\lambda \vDash l} (-1)^{l_\lambda} B_{\lambda_1}(x, t_1, ..., t_{\lambda_1}) \prod_{j=2}^{l_\lambda} A_{\lambda_j}(x, t_{|\lambda|_{j-1}+1}, ..., t_{|\lambda|_j}),$$

which is precisely formula (3.4). □

Notice that $y_1$ is given directly by the solution of the ODE $\dot{y}_1(x, t_1) = a_1(x) y(x, t_1) + a_{0,1}(x, t_1)$.

In this paper we will always consider the case $u(0, t, s) = t_1\, s$. In the Riccati-type equation (3.2) we have $a_2(x) = i x^{-2}$, $a_1 = 0$, $a_0(x, t, s) = -i x^2\, e^{i x t_1}\, s$. Therefore the corresponding second order linear ODE writes as

$$\ddot{U} + \frac{2}{x} \dot{U} + (e^{i x t_1} s) * U = 0. \tag{3.6}$$

with $U(0, t, s) = t_1\, s$ and $\dot{U}(0, t, s) = 0$. Then according to formula (3.4) we infer the expression

$$v_l(x, t_1, ..., t_l) = -i x^2 \sum_{\lambda \vDash l} (-1)^{l_\lambda} \frac{\dot{U}_{\lambda_1+1}(x, t_1, ..., t_{\lambda_1+1})}{t_{\lambda_1+1}} \prod_{j=2}^{l_\lambda} \frac{U_{\lambda_j+1}(x, t_{|\lambda|_{j-1}+1}, ..., t_{|\lambda|_j+1})}{t_{|\lambda|_j+1}}, \tag{3.7}$$

for all $l \geqslant 1$, for the solution $v$ of the equation (3.2).

In order to compute the solution $U$ we need an other elementary lemma.

**Lemma 3.8.** *In the set up of the previous lemma, any ODE*

$$\ddot{u} + a_1(x)\dot{u} + a_0 * u \;=\; f(x,t,s)\,,$$

*is equivalent to the ODE*

$$\ddot{w} - \frac{1}{4}(2\,\dot{a}_1 + a_1^2 - 4\,a_0) * w \;=\; e^{\frac{1}{2}\int a_1} f,$$

*via the identification*

$$u \;=\; e^{-\frac{1}{2}\int a_1}\, w.$$

**Proof.** We write $u = A(x)w$ and we observe the elementary equalities $\dot{u} = \dot{A}\,w + A\,\dot{w}$,

$$\ddot{u} \;=\; \ddot{A}\,w + 2\,\dot{A}\,\dot{w} + A\,\ddot{w}\,.$$

Then the ODE on $u$ in the statement is equivalent to the ODE

$$A\,\ddot{w} + 2\,\dot{A}\,\dot{w} + \ddot{A}\,w + a_1 A\,\dot{w} + a_1 \dot{A}\,w + A\,a_0 * w \;=\; f\,,$$

i.e.

$$A\,\ddot{w} + (2\,\dot{A} + a_1 A)\dot{w} + (\ddot{A} + a_1\dot{A} + A\,a_0) * w \;=\; f\,.$$

We seek for $A$ solution of the ODE $2\,\dot{A} + a_1 A = 0$, i.e. $A = e^{-\frac{1}{2}\int a_1}$. We infer in particular the identity

$$\ddot{A} + \frac{1}{2}\,\dot{a}_1 A + \frac{1}{2}\,a_1 \dot{A} \;=\; 0\,.$$

Using this we write

$$\begin{aligned} \ddot{A} + a_1\dot{A} + A\,a_0 \;&=\; -\frac{1}{2}\,\dot{a}_1 A + \frac{1}{2}\,a_1\dot{A} + A\,a_0 \\ &=\; -\frac{A}{4}\,(2\,\dot{a}_1 + a_1^2 - 4\,a_0)\,, \end{aligned}$$

which implies the required conclusion. □

We infer that if we set $U = x^{-1} w$, then the equation (3.6) is equivalent to the equation

$$\ddot{w} + (e^{ixt_1}\,s) * w = 0, \tag{3.8}$$

with the initial conditions $w(0,t,s)\equiv 0$ and $\dot{w}(0,t,s)=u(0,t,s)=t_1 s$. (The later condition follows deriving the equality $w=xU$). Then formula (3.7) implies the expression

$$v_l(x,t_1,...,t_l) \;=\; -ix\sum_{\lambda\vDash l}(-1)^{l_\lambda}\frac{\dot{w}_{\lambda_1+1}(x,t_1,...,t_{\lambda_1+1})}{t_{\lambda_1+1}}\prod_{j=2}^{l_\lambda}\frac{w_{\lambda_j+1}(x,t_{|\lambda|_{j-1}+1},...,t_{|\lambda|_j+1})}{x\,t_{|\lambda|_j+1}}$$

$$+\; ix\sum_{\lambda\vDash l}(-1)^{l_\lambda}\prod_{j=1}^{l_\lambda}\frac{w_{\lambda_j+1}(x,t_{|\lambda|_{j-1}+1},...,t_{|\lambda|_j+1})}{x\,t_{|\lambda|_j+1}},$$

for all $l\geqslant 1$, for the solution $v$ of the equation (3.2).

## 3.2 Expression for the solution of the equation (3.8)

We notice first the equalities $w(0,t,s)=\ddot{w}(0,t,s)=0$, $\dot{w}_1(0,t_1)=t_1$, $\ddot{w}_1\equiv 0$. We set

$$w_{p,l}(t_1,...,t_l) \;:=\; \frac{1}{p!}\,w_l^{(p)}(0,t_1,...,t_l)\,.$$

Deriving the equation (3.8) with respect to the variable $x$ and evaluating at $x=0$ we infer for all $l\geqslant 2$

$$w_{k+2,l}(t_1,...,t_l) \;=\; -\frac{1}{(k+2)(k+1)}\sum_{r=0}^{k-1}\frac{i^r}{r!}\,t_1^r\,w_{k-r,l-1}(t_2,...,t_l)\,,$$

which we rewrite as

$$\begin{aligned} w_{p,l}(t_1,...,t_l) \;&=\; -\frac{1}{p(p-1)}\sum_{r=0}^{p-3}\frac{i^r}{r!}\,t_1^r\,w_{p-r-2,l-1}(t_2,...,t_l)\\ &=\; \frac{1}{p(p-1)}\sum_{r=2}^{p-1}\frac{i^r}{(r-2)!}\,t_1^{r-2}\,w_{p-r,l-1}(t_2,...,t_l)\,. \end{aligned}$$

We deduce

$$w_{p,2}(t_1,t_2) \;=\; i^{p-1}\frac{p-2}{p!}\,t_1^{p-3}\,t_2\,,$$

for all $p\geqslant 3$ and zero otherwise. We deduce also

$$\begin{aligned} w_{p,3}(t_1,t_2,t_3) \;&=\; \frac{1}{p(p-1)}\sum_{r=2}^{p-1}\frac{i^r}{(r-2)!}\,t_1^{r-2}\,w_{p-r,2}(t_2,t_3)\\ &=\; \frac{1}{p(p-1)}\sum_{r=2}^{p-3}\frac{i^r}{(r-2)!}\,t_1^{r-2}\,w_{p-r,2}(t_2,t_3)\,, \end{aligned}$$

since $w_{p-r,2}=0$ for $p-r\geqslant 3$. Thus

$$w_{p,3}(t_1,t_2,t_3) \;=\; \frac{i^{p-1}}{p(p-1)}\sum_{r=2}^{p-3}(p-r-2)\frac{t_1^{r-2}\,t_2^{p-r-3}}{(r-2)!(p-r)!}\,t_3\,.$$

We compute now

$$w_{p,4}(t_1,...,t_4) = \frac{1}{p(p-1)}\sum_{r=2}^{p-1}\frac{i^r}{(r-2)!}\,t_1^{r-2}\,w_{p-r,3}(t_2,...,t_4)$$

$$= \frac{1}{p(p-1)}\sum_{r=2}^{p-5}\frac{i^r}{(r-2)!}\,t_1^{r-2}\,w_{p-r,3}(t_2\,,...,t_4)\,,$$

since $w_{p-r,3}=0$ for $p-r\geqslant 5$. Thus

$$w_{p,4}(t_1,...,t_4) = \frac{i^{p-1}}{p(p-1)}\sum_{r_1=2}^{p-5}\sum_{r_2=2}^{p-r_1-3}\frac{(p-r_1-r_2-2)\,t_1^{r_1-2}\,t_2^{r_2-2}\,t_3^{p-r_1-r_2-3}}{(p-r_1)(p-r_1-1)(r_1-2)!(r_2-2)!(p-r_1-r_2)!}\,t_4\,.$$

We compute now

$$w_{p,5}(t_1,...,t_5) = \frac{1}{p(p-1)}\sum_{r=2}^{p-1}\frac{i^r}{(r-2)!}\,t_1^{r-2}\,w_{p-r,4}(t_2,...,t_5)$$

$$= \frac{1}{p(p-1)}\sum_{r=2}^{p-7}\frac{i^r}{(r-2)!}\,t_1^{r-2}\,w_{p-r,4}(t_2\,,...,t_5)\,,$$

since $w_{p-r,4}=0$ for $p-r\geqslant 7$. Thus

$$w_{p,5}(t_1,...,t_5)$$

$$= \frac{i^{p-1}}{p(p-1)}\sum_{r_1=2}^{p-7}\sum_{r_2=2}^{p-r_1-5}\sum_{r_3=2}^{p-r_1-r_2-3}\frac{p-r_1-r_2-r_3-2}{(p-r_1)(p-r_1-1)(p-r_1-r_2)(p-r_1-r_2-1)}\times$$

$$\times\ \frac{t_1^{r_1-2}\,t_2^{r_2-2}\,t_3^{r_3-2}\,t_4^{p-r_1-r_2-r_3-3}}{(r_1-2)!(r_2-2)!(r_3-2)!(p-r_1-r_2-r_3)!}\,t_5\,.$$

We infer the general formula

$$w_{p,l}(t_1,...,t_l) = \frac{i^{p-1}}{p(p-1)}\sum_{\substack{r_j=2\\1\leqslant j\leqslant l-2}}^{p-|r|_{j-1}-2(l-j)+1}(p-|r|_{l-2}-2)\times$$

$$\times\ \left[\prod_{j=1}^{l-3}\frac{1}{(p-|r|_j)(p-|r|_j-1)}\right]\left[\prod_{j=1}^{l-2}\frac{t_j^{r_j-2}}{(r_j-2)!}\right]\frac{t_{l-1}^{p-|r|_{l-2}-3}}{(p-|r|_{l-2})!}\,t_l\,,$$

for all $l\geqslant 3$. Performing the change of variables $r_j'=r_j-2$ we infer

$$w_{p,l}(t_1,...,t_l) = \frac{i^{p-1}}{p(p-1)}\sum_{\substack{r_j=0\\1\leqslant j\leqslant l-2}}^{p-|r|_{j-1}-2l+1}(p-|r|_{l-2}-2l+2)\times$$

$$\times\ \left[\prod_{j=1}^{l-3}\frac{1}{(p-|r|_j-2j)(p-|r|_j-2j-1)}\right]\left[\prod_{j=1}^{l-2}\frac{t_j^{r_j}}{r_j!}\right]\frac{t_{l-1}^{p-|r|_{l-2}-2l+1}}{(p-|r|_{l-2}-2l+4)!}\,t_l\,.$$

If we set $r_{l-1} := p - |r|_{l-2} - 2l + 1$, then we can rewrite the previous sum as

$$w_{p,l}(t_1, ..., t_l) = \frac{i^{p-1}}{p(p-1)} \sum_{r \vDash_0^{l-1} p-2l+1} (r_{l-1}+1) \times$$

$$\times \left[\prod_{j=1}^{l-3} \frac{1}{(p-|r|_j-2j)(p-|r|_j-2j-1)}\right]\left[\prod_{j=1}^{l-2} \frac{t_j^{r_j}}{r_j!}\right]\frac{t_{l-1}^{r_{l-1}}}{(r_{l-1}+3)!} t_l$$

$$= \frac{i^{p-1}}{p(p-1)} \sum_{r \vDash_0^{l-1} p-2l+1} \left[\prod_{j=1}^{l-2} \frac{1}{(p-|r|_j-2j)(p-|r|_j-2j-1)}\right]\left[\prod_{j=1}^{l-1} \frac{t_j^{r_j}}{r_j!}\right] t_l\,,$$

where $r \vDash_0^{l-1} p - 2l + 1$ denotes the compositions of length $l-1$ and $|r| = p - 2l + 1$.

We notice that the previous formula hold also for $l = 2$. In conclusion

$$w_1(x, t_1) = x t_1\,,$$

$$w_l(x, t_1, ..., t_l) = \sum_{p \geqslant 2l-1} w_{p,l}(t_1, ..., t_l)\, x^p\,,$$

for all $l \geqslant 2$.

### 3.3 The explicit expression of the coefficients $C_D$

By the result of the previous subsection, we infer the expression

$$\frac{w_{\lambda_j+1}(x, t_{|\lambda|_{j-1}+1}, ..., t_{|\lambda|_j+1})}{t_{|\lambda|_j+1}} = \sum_{p \geqslant 2\lambda_j+1} \frac{w_{p,\lambda_j+1}(t_{|\lambda|_{j-1}+1}, ..., t_{|\lambda|_j+1})}{t_{|\lambda|_j+1}} x^p\,,$$

with

$$\frac{w_{p,\lambda_j+1}(t_{|\lambda|_{j-1}+1}, ..., t_{|\lambda|_j+1})}{t_{|\lambda|_j+1}}$$

$$= \frac{i^{p-1}}{p(p-1)} \sum_{r \vDash_0^{\lambda_j} p-2\lambda_j-1} \left[\prod_{s=1}^{\lambda_j-1} \frac{1}{(p-|r|_s-2s)(p-|r|_s-2s-1)}\right]\left[\prod_{s=1}^{\lambda_j} \frac{t_{|\lambda|_{j-1}+s}^{r_s}}{r_s!}\right].$$

Thus

$$\frac{w_{\lambda_j+1}(x, t_{|\lambda|_{j-1}+1}, ..., t_{|\lambda|_j+1})}{t_{|\lambda|_j+1}}$$

$$= \sum_{p \geqslant 2\lambda_j+1} \frac{i^{p-1}x^p}{p(p-1)} \sum_{r \vDash_0^{\lambda_j} p-2\lambda_j-1} \left[\prod_{s=1}^{\lambda_j-1} \frac{1}{(p-|r|_s-2s)(p-|r|_s-2s-1)}\right]\left[\prod_{s=1}^{\lambda_j} \frac{t_{|\lambda|_{j-1}+s}^{r_s}}{r_s!}\right],$$

and

$$\frac{\dot{w}_{\lambda_1+1}(x,t_1,...,t_{\lambda_1+1})}{t_{\lambda_1+1}}$$

$$= \frac{1}{x}\sum_{p\geqslant 2\lambda_1+1}\frac{i^{p-1}x^p}{(p-1)}\sum_{r\vDash_0^{\lambda_1} p-2\lambda_1-1}\left[\prod_{s=1}^{\lambda_1-1}\frac{1}{(p-|r|_s-2s)(p-|r|_s-2s-1)}\right]\left[\prod_{s=1}^{\lambda_1}\frac{t_s^{r_s}}{r_s!}\right].$$

We infer in conclusion

$$v_l(x,t_1,...,t_l)$$

$$= -i\sum_{\lambda\vDash l}\frac{(-1)^{l_\lambda}}{x^{l_\lambda-1}}\sum_{\substack{p_j\geqslant 2\lambda_j+1\\ r_j\vDash_0^{\lambda_j} p_j-2\lambda_j-1\\ j=1,...l_\lambda}}(p_1-1)\prod_{j=1}^{l_\lambda}\frac{i^{p_j-1}x^{p_j}}{p_j(p_j-1)}\times$$

$$\times\left[\prod_{j=1}^{l_\lambda}\prod_{s_j=1}^{\lambda_j-1}\frac{1}{(p_j-|r_j|_{s_j}-2s_j)(p_j-|r_j|_{s_j}-2s_j-1)}\right]\left[\prod_{j=1}^{l_\lambda}\prod_{s_j=1}^{\lambda_j}\frac{t_{|\lambda|_{j-1}+s_j}^{r_{j,s_j}}}{r_{j,s_j}!}\right],$$

for all $l\geqslant 1$. Thus

$$v_l(x,t_1,...,t_l)$$

$$= (-1)^{l+1}i\sum_{H\in\mathbb{Z}^l_{\geqslant 0}}\sum_{\lambda\vDash l}(-1)^{l_\lambda}(|H|_{\lambda_1}+2\lambda_1)\,i^{|H|}\,x^{|H|+2l+1}\times$$

$$\times\left[\prod_{j=1}^{l_\lambda}\prod_{s=|\lambda|_{j-1}}^{|\lambda|_j-1}\frac{1}{[\sum_{r=s+1}^{|\lambda|_j}h_r+2(|\lambda|_j-s)+1][\sum_{r=s+1}^{|\lambda|_j}h_r+2(|\lambda|_j-s)]}\right]\frac{t^H}{H!},$$

for all $l\geqslant 1$, which we rewrite as

$$v_l(x,t_1,...,t_l)$$

$$= (-1)^{l+1}i\sum_{H\in\mathbb{Z}^l_{\geqslant 0}}i^{|H|}\,x^{|H|+2l+1}\frac{t^H}{H!}\sum_{\lambda\vDash l}(-1)^{l_\lambda}(|H|_{\lambda_1}+2\lambda_1)\times$$

$$\times\left[\prod_{j=1}^{l_\lambda}\prod_{s=|\lambda|_{j-1}}^{|\lambda|_j-1}\frac{1}{[\sum_{r=s+1}^{|\lambda|_j}h_r+2(|\lambda|_j-s)+1][\sum_{r=s+1}^{|\lambda|_j}h_r+2(|\lambda|_j-s)]}\right].$$

We remind now that thanks to lemma 3.6 the solution $u$ of the Riccati-type equation in lemma 3.4 is given by

$$u_l(x,t_1,...,t_l) = -i\,e^{-ix\sum_{j=1}^l t_j}v_l(x,t_1,...,t_l).$$

Then using the expansion

$$e^{-ix\sum_{j=1}^l t_j} = \sum_{H\in\mathbb{Z}^l_{\geqslant 0}}(-i)^{|H|}x^{|H|}\frac{t^H}{H!},$$

we infer the formula

$$
\begin{aligned}
& u_l(x, t_1, ..., t_l) \\
&= (-1)^{l+1} \sum_{D\in\mathbb{Z}^l_{\geqslant 0}} (-i)^{|D|} x^{|D|+2l+1} t^D \sum_{0\leqslant H\leqslant D} \frac{(-1)^{|H|}}{H!\,(D-H)!} \sum_{\lambda\vDash l} (-1)^{l_\lambda} (|H|_{\lambda_1} + 2\lambda_1) \times \\
&\times \left[ \prod_{j=1}^{l_\lambda} \prod_{s=|\lambda|_{j-1}}^{|\lambda|_j - 1} \frac{1}{[\sum_{r=s+1}^{|\lambda|_j} h_r + 2(|\lambda|_j - s) + 1][\sum_{r=s+1}^{|\lambda|_j} h_r + 2(|\lambda|_j - s)]} \right],
\end{aligned}
$$

for all $l \geqslant 1$. We conclude the explicit expression

$$
\begin{aligned}
C_D &= (-1)^{l+1} (-i)^{|D|} (|D| + 2l + 1)! \sum_{0\leqslant H\leqslant D} \frac{(-1)^{|H|}}{H!\,(D-H)!} \sum_{\lambda\vDash l} (-1)^{l_\lambda} (|H|_{\lambda_1} + 2\lambda_1) \times \\
&\times \left[ \prod_{j=1}^{l_\lambda} \prod_{s=|\lambda|_{j-1}}^{|\lambda|_j - 1} \frac{1}{[\sum_{r=s+1}^{|\lambda|_j} h_r + 2(|\lambda|_j - s) + 1][\sum_{r=s+1}^{|\lambda|_j} h_r + 2(|\lambda|_j - s)]} \right],
\end{aligned}
$$

for all $D \in \mathbb{Z}^l_{\geqslant 0}$. This is precisely the expression of the coefficients $C_D$ in the statement of theorem 1.6.

# 4 Expliciting the integrability equations

In this section we provide some basic tools for a general expression of the integrability equations

$$
\operatorname{Circ}\operatorname{Sym}_{3,\ldots,k+1} \Theta_k^{\nabla} = 0,
$$

for $k \geqslant 4$. We notice in particular that the equation $\operatorname{Circ}\operatorname{Sym}_{3,4,5}\Theta_4(g) = 0$ writes as

$$
\operatorname{Circ}\operatorname{Sym}_{3,4,5}\big[3 d_1^{\nabla^g} (\nabla^g R^g)_2 - 2\,\tilde{R}^g \wedge_1 \tilde{R}^g\,\big] = 0,
$$

with $\tilde{R}^g := \operatorname{Sym}_{2,3} R^g$. Its vanishing has been proved in [Pali] by using a direct computation.

## 4.1 Expliciting the powers of the 1-differential. Part II

From now on let $\rho_{k,j} := \big(R^{\nabla}.R^{(k+1-j)}\big)^{(j-2)}$ and notice that

$$
\rho_{k,j} \in C^\infty\big(M, T_M^{*,\otimes(j-2)} \otimes_{\mathbb{R}} \Lambda^2 T_M^* \otimes_{\mathbb{R}} T_M^{*,\otimes(k+1-j)} \otimes_{\mathbb{R}} \Lambda^2 T_M^* \otimes_{\mathbb{R}} T_M^* \otimes_{\mathbb{R}} \mathbb{C}T_M\big),
$$

satisfies the circular identity with respect to its last three entries and in the case $j \leqslant k$ the tensor $\rho_{k,j} \cdot v$ satisfies the circular identity with respect to its last three entries for any vector $v \in T_M$. We remind the identity

$$
\begin{aligned}
R^{(k+1)}(2, ..., k+3, 1, k+4) &= R^{(k+1)}(3, ..., k+2, 2, k+3, 1, k+4) \\
&+ \sum_{j=2}^{k+1} \rho_{k,j}(3, ..., j, 2, j+1, ..., k+3, 1, k+4),
\end{aligned}
$$

obtained in the proof of the identity (2.3) in proposition 2.5. We rewrite it in the more general form

$$R^{(c-b+1)}(1,b,\ldots,c,2,d) \;=\; R^{(c-b+1)}(b,\ldots,c-1,1,c,2,d)$$

$$+\sum_{h=b-1}^{c-2}\rho_{c-b,h-b+3}(b,\ldots,h,1,h+1,\ldots,c,2,d)\,,$$

for any $3\leqslant b<c$. Deriving we obtain the general commutation formula

$$R^{(c-a+1)}(a,\ldots,b-1,1,b,\ldots,c,2,d) \;=\; R^{(c-a+1)}(a,\ldots,c-1,1,c,2,d)$$

$$+\sum_{h=b-1}^{c-2}\rho_{c-a,h-a+3}(a,\ldots,h,1,h+1,\ldots,c,2,d)\,, \tag{4.1}$$

for all $3\leqslant a\leqslant b<c<d$. We show now two key elementary lemmas.

**Lemma 4.1.** *In the set up of proposition 2.5, the identity hold*

$$\mathrm{Circ}\,\mathrm{Sym}_{3,\ldots,k+4}T_k \;=\; \mathrm{Circ}\,\mathrm{Sym}_{4,\ldots,k+4}\hat{T}_k\,,$$

*with*

$$\hat{T}_k(1,\ldots,k+4) \;:=\; (-1)^k\sum_{j=2}^{k+1}(j-1)\,\rho_{k,j}(4,\ldots,j+1,3,j+2,\ldots,k+3,1,2,k+4)\,.$$

**Proof.** We notice first that for any tensor $\theta\equiv\theta(1,\ldots,k+4)$, the tensor

$$\tilde{\theta}(1,\ldots,k+4) \;:=\; \sum_{j=3}^{k+4}\theta(1,2,4,\ldots,j,3,j+1,\ldots,k+4)\,,$$

satisfies the identity

$$\mathrm{Sym}_{3,\ldots,k+4}\,\theta \;=\; \mathrm{Sym}_{4,\ldots,k+4}\,\tilde{\theta}\,.$$

Applying the previous definition to our tensor $T_k$ we obtain the expression

$$\tilde{T}_k(1,\ldots,k+4) \;=\; (-1)^k\sum_{j=3}^{k+3}R^{(k+1)}(4,\ldots,j,3,j+1,\ldots,k+3,1,2,k+4)$$

$$+\;(-1)^k\,R^{(k+1)}(4,\ldots,k+4,1,2,3)\,.$$

Using the commutation identity (4.1) we infer

$$R^{(k+1)}(4,\ldots,j,3,j+1,\ldots,k+3,1,2,k+4)$$

$$=\;R^{(k+1)}(4,\ldots,k+3,3,1,2,k+4)+\sum_{h=j}^{k+2}\rho_{k,h-3}(4,\ldots,h,3,h+1,\ldots,k+3,1,2,k+4)\,.$$

We infer the equality

$$\begin{aligned}
&\tilde{T}_k(1, ..., k+4)\\
= &\ (-1)^k (k+1) R^{(k+1)}(4, ..., k+3, 3, 1, 2, k+4) + (-1)^k R^{(k+1)}(4, ..., k+4, 1, 2, 3)\\
+ &\ (-1)^k \sum_{j=3}^{k+2} \sum_{h=j}^{k+2} \rho_{k,h-3}(4, ..., h, 3, h+1, ..., k+3, 1, 2, k+4)\,.
\end{aligned}$$

Using the elementary equality

$$\sum_{s=a}^{b} \sum_{h=s}^{b} C_h = \sum_{h=a}^{b} (h-a+1)\, C_h\,,$$

we obtain the identity

$$\begin{aligned}
&\tilde{T}_k(1, ..., k+4)\\
= &\ (-1)^k (k+1) R^{(k+1)}(4, ..., k+3, 3, 1, 2, k+4) + (-1)^k R^{(k+1)}(4, ..., k+4, 1, 2, 3)\\
+ &\ (-1)^k \sum_{h=3}^{k+2} (h-2)\, \rho_{k,h-3}(4, ..., h, 3, h+1, ..., k+3, 1, 2, k+4)\,.
\end{aligned}$$

Using the algebraic and differential Bianchi identities and performing the change of variable $h = j+1$, we infer the required conclusion. □

**Lemma 4.2.** *In the set up of proposition 2.5, for all integers $k \geqslant 1$ hold the identity*

$$\operatorname{Circ}\operatorname{Sym}_{3,...,k+4}\Big[(d_1^{\nabla})^k (\nabla R^{\nabla})_2\Big] = \sum_{j=2}^{k+1} \operatorname{Circ}\operatorname{Sym}_{4,...,k+4} (\check{T}_{k,j} + \check{Q}_{k,j} + \check{V}_{k,j})\,,$$

*with*

$$\begin{aligned}
\check{T}_{k,j}(1, ..., k+4) &:= (-1)^k (j-2)\, \rho_{k,j}(4, ..., j+1, 3, j+2, ..., k+3, 1, 2, k+4)\,,\\
\check{Q}_{k,j}(1, ..., k+4) &:= -(-1)^k \sum_{s=3}^{j} \rho_{k,j}(4, ..., s, 3, s+1, ..., j, 2, j+1, ..., k+3, 1, k+4)\\
&\quad - (-1)^k \sum_{s=j+1}^{k+3} \rho_{k,j}(4, ..., j+1, 2, j+2, ..., s, 3, s+1, ..., k+3, 1, k+4)\,,\\
\check{V}_{k,j}(1, ..., k+4) &:= (-1)^k \sum_{s=3}^{j} \rho_{k,j}(4, ..., s, 3, s+1, ..., j, 1, j+1, ..., k+3, 2, k+4)\\
&\quad + (-1)^k \sum_{s=j+1}^{k+3} \rho_{k,j}(4, ..., j+1, 1, j+2, ..., s, 3, s+1, ..., k+3, 2, k+4)\,.
\end{aligned}$$

**Proof.** In the computation that will follow we will denote by $\rho_{k,j}(...)_p$ with $p=0,1$, the terms that summed-up together annihilate the operator Circ $\mathrm{Sym}_{4,...,k+4}$. For all $2 \leqslant j \leqslant k+1$, we define

$$\begin{aligned}
\hat{T}_{k,j}(1,...,k+4) &:= (-1)^k (j-1)\, \rho_{k,j}(4,...,j+1,3,j+2,...,k+3,1,2,k+4)\,,\\
\hat{Q}_{k,j}(1,...,k+4) &:= -(-1)^k \sum_{s=3}^{j} \rho_{k,j}(4,...,s,3,s+1,...,j,2,j+1,...,k+3,1,k+4)\\
&- (-1)^k \sum_{s=j+1}^{k+3} \rho_{k,j}(4,...,j+1,2,j+2,...,s,3,s+1,...,k+3,1,k+4)\\
&- (-1)^k \rho_{k,j}(4,...,j+1,2,j+2,...,k+4,1,3)_2\,,\\
\hat{V}_{k,j}(1,...,k+4) &:= (-1)^k \sum_{s=3}^{j} \rho_{k,j}(4,...,s,3,s+1,...,j,1,j+1,...,k+3,2,k+4)\\
&+ (-1)^k \sum_{s=j+1}^{k+3} \rho_{k,j}(4,...,j+1,1,j+2,...,s,3,s+1,...,k+3,2,k+4)\\
&= (-1)^k \rho_{k,j}(4,...,j+1,1,j+2,...,k+4,2,3)_1\,,
\end{aligned}$$

and we notice the equalities

$$\begin{aligned}
\mathrm{Sym}_{3,...,k+4}\, Q_k &= \sum_{j=2}^{k+1} \mathrm{Sym}_{4,...,k+4}\, \hat{Q}_{k,j}\,,\\
\mathrm{Sym}_{3,...,k+4}\, V_k &= \sum_{j=2}^{k+1} \mathrm{Sym}_{4,...,k+4}\, \hat{V}_{k,j}\,.
\end{aligned}$$

Then proposition 2.5 and lemma 4.1 imply

$$\mathrm{Circ}\,\mathrm{Sym}_{3,...,k+4}\Big[(d_1^{\nabla})^k\,(\nabla R^{\nabla})_2\Big] = \sum_{j=2}^{k+1} \mathrm{Circ}\,\mathrm{Sym}_{4,...,k+4}\,(\hat{T}_{k,j} + \hat{Q}_{k,j} + \hat{V}_{k,j})\,. \tag{4.2}$$

We write

$$\begin{aligned}
\hat{T}_{k,j}(1,...,k+4) &= (-1)^k (j-2)\, \rho_{k,j}(4,...,j+1,3,j+2,...,k+3,1,2,k+4)\\
&- (-1)^k \rho_{k,j}(4,...,j+1,3,j+2,...,k+4,1,2)\\
&- (-1)^k \rho_{k,j}(4,...,j+1,3,j+2,...,k+3,2,k+4,1)\\
&= (-1)^k (j-2)\, \rho_{k,j}(4,...,j+1,3,j+2,...,k+3,1,2,k+4)\\
&- (-1)^k \rho_{k,j}(4,...,j+1,3,j+2,...,k+4,1,2)_1\\
&+ (-1)^k \rho_{k,j}(4,...,j+1,3,j+2,...,k+4,2,1)_2\,.
\end{aligned}$$

We infer the identity

$$
\begin{aligned}
&\operatorname{Circ}\operatorname{Sym}_{4,\dots,k+4}\left(\hat{T}_{k,j}+\hat{Q}_{k,j}+\hat{V}_{k,j}\right)\\
=\;&\operatorname{Circ}\operatorname{Sym}_{4,\dots,k+4}\left(\check{T}_{k,j}+\check{Q}_{k,j}+\check{V}_{k,j}\right),
\end{aligned}
$$

which combined with the identity (4.2) implies the required statement. □

**Lemma 4.3.** *In the set up of proposition 2.5, for all integers $k \geqslant 1$ hold the identity*

$$
\operatorname{Circ}\operatorname{Sym}_{3,\dots,k+4}\Big[(d_1^{\nabla})^k\,(\nabla R^{\nabla})_2\Big] \;=\; (-1)^k \sum_{j=2}^{k+1} \operatorname{Circ}\operatorname{Sym}_{4,\dots,k+4}\big[(j-2)\,\mathbb{T}_{k,j} - \operatorname{Alt}_2 \mathbb{Q}_{k,j}\big],
$$

*with*

$$
\begin{aligned}
&\mathbb{T}_{k,j}(1,\dots,k+4)\\
:=\;& \sum_{I\subseteq\{4,\dots,j+1\}} R^{(|I|)}(I,3,j+2,R^{(k+1-j+|\complement I|)}(\complement I,j+3,\dots,k+3,1,2,k+4))\\
-\;& \sum_{\substack{I\subseteq\{4,\dots,j+1\}\\ j+3\leqslant h\leqslant k+3}} R^{(k+1-j+|I|)}(I,j+3,\dots,R^{(|\complement I|)}(\complement I,3,j+2,h),\dots,k+3,1,2,k+4)\\
-\;& \sum_{I\subseteq\{4,\dots,j+1\}} R^{(k+1-j+|I|)}(I,j+3,\dots,k+3,R^{(|\complement I|)}(\complement I,3,j+2,1),2,k+4)\\
-\;& \sum_{I\subseteq\{4,\dots,j+1\}} R^{(k+1-j+|I|)}(I,j+3,\dots,k+3,1,R^{(|\complement I|)}(\complement I,3,j+2,2),k+4)\\
-\;& \sum_{I\subseteq\{4,\dots,j+1\}} R^{(k+1-j+|I|)}(I,j+3,\dots,k+3,1,2,R^{(|\complement I|)}(\complement I,3,j+2,k+4)),
\end{aligned}
$$

*and with*

$$
\begin{aligned}
&\mathbb{Q}_{k,j}(1,\dots,k+4)\\
:=\;& \sum_{\substack{3\leqslant s\leqslant j\\ I\subseteq\{4,\dots,s,3,s+1,\dots,j\}}} R^{(|I|)}(I,2,j+1,R^{(k+1-j+|\complement I|)}(\complement I,j+2,\dots,k+3,1,k+4))\\
-\;& \sum_{\substack{3\leqslant s\leqslant j\\ I\subseteq\{4,\dots,s,3,s+1,\dots,j\}\\ j+2\leqslant h\leqslant k+3}} R^{(k+1-j+|I|)}(I,j+2,\dots,R^{(|\complement I|)}(\complement I,2,j+1,h),\dots,k+3,1,k+4)\\
-\;& \sum_{\substack{3\leqslant s\leqslant j\\ I\subseteq\{4,\dots,s,3,s+1,\dots,j\}}} R^{(k+1-j+|I|)}(I,j+2,\dots,k+3,R^{(|\complement I|)}(\complement I,2,j+1,1),k+4)\\
-\;& \sum_{\substack{3\leqslant s\leqslant j\\ I\subseteq\{4,\dots,s,3,s+1,\dots,j\}}} R^{(k+1-j+|I|)}(I,j+2,\dots,k+3,1,R^{(|\complement I|)}(\complement I,2,j+1,k+4))
\end{aligned}
$$

$$+ \sum_{I \subseteq \{4,...,j+1\}} R^{(|I|)}(I, 2, 3, R^{(k+1-j+|\complement I|)}(\complement I, j+2, ..., k+3, 1, k+4))$$

$$- \sum_{\substack{I \subseteq \{4,...,j+1\} \\ j+2 \leqslant h \leqslant k+3}} R^{(k+1-j+|I|)}(I, j+2, ..., R^{(|\complement I|)}(\complement I, 2, 3, h), ..., k+3, 1, k+4)$$

$$- \sum_{I \subseteq \{4,...,j+1\}} R^{(k+1-j+|I|)}(I, j+2, ..., k+3, 1, R^{(|\complement I|)}(\complement I, 2, 3, k+4))$$

$$+ \sum_{\substack{j+2 \leqslant s \leqslant k+3 \\ I \subseteq \{4,...,j+1\}}} R^{(|I|)}(I, 2, j+2, R^{(k+1-j+|\complement I|)}(\complement I, j+3, ..., s, 3, s+1, ..., k+3, 1, k+4))$$

$$- \sum_{\substack{j+3 \leqslant s \leqslant k+3 \\ I \subseteq \{4,...,j+1\} \\ j+3 \leqslant h \leqslant s}} R^{(k+1-j+|I|)}(I, j+3, ..., R^{(|\complement I|)}(\complement I, 2, j+2, h), ..., s, 3, s+1, ..., k+3, 1, k+4)$$

$$- \sum_{\substack{j+2 \leqslant s \leqslant k+3 \\ I \subseteq \{4,...,j+1\}}} R^{(k+1-j+|I|)}(I, j+3, ..., s, R^{(|\complement I|)}(\complement I, 2, j+2, 3), s+1, ..., k+3, 1, k+4)$$

$$- \sum_{\substack{j+2 \leqslant s \leqslant k+2 \\ I \subseteq \{4,...,j+1\} \\ s+1 \leqslant h \leqslant k+3}} R^{(k+1-j+|I|)}(I, j+3, ..., s, 3, s+1, ..., R^{(|\complement I|)}(\complement I, 2, j+2, h), ..., k+3, 1, k+4)$$

$$- \sum_{\substack{j+2 \leqslant s \leqslant k+3 \\ I \subseteq \{4,...,j+1\}}} R^{(k+1-j+|I|)}(I, j+3, ..., s, 3, s+1, ..., k+3, R^{(|\complement I|)}(\complement I, 2, j+2, 1), k+4)$$

$$- \sum_{\substack{j+2 \leqslant s \leqslant k+3 \\ I \subseteq \{4,...,j+1\}}} R^{(k+1-j+|I|)}(I, j+3, ..., s, 3, s+1, ..., k+3, 1, R^{(|\complement I|)}(\complement I, 2, j+2, k+4))\,,$$

*where the symbol* $I \subseteq \{4, ..., s, 3, s+1, ..., j\}$ *means a subset of the elements* $4, ..., s, 3, s+1, ..., j$ *written in this order and* $\complement I \subseteq \{4, ..., s, 3, s+1, ..., j\}$ *denotes its complementary set respecting the same order.*

**Proof.** We remind first the definition

$$(R.R^{(p)})(1, ..., p+5) = R(1, 2, R^{(p)}(3, ..., p+5)) - \sum_{h=3}^{p+5} R^{(p)}(3, ..., R(1, 2, h), ..., p+5)\,,$$

which writes in the case

$$(R.R^{(k+1-j)})(j-1, ..., k+4) = R(j-1, j, R^{(k+1-j)}(j+1, ..., k+4)) - \sum_{h=j+1}^{k+4} R^{(k+1-j)}(j+1, ..., R(j-1, j, h), ..., k+4)\,.$$

We infer the expression

$$\begin{aligned}
&\rho_{k,j}(1,...,k+4)\\
&= \sum_{I\subseteq\{1,...,j-2\}} R^{(|I|)}(I,j-1,j,R^{(k+1-j+|\complement I|)}(\complement I,j+1,...,k+4))\\
&- \sum_{I\subseteq\{1,...,j-2\}} \sum_{h=j+1}^{k+4} R^{(k+1-j+|I|)}(I,j+1,...,R^{(|\complement I|)}(\complement I,j-1,j,h),...,k+4)\,.
\end{aligned}$$

Performing the change of variables

$$\begin{matrix}
1 & ... & j-2 & j-1 & j & ... & k+1 & k+2 & k+3 & k+4\\
\downarrow & & \downarrow & \downarrow & \downarrow & & \downarrow & \downarrow & \downarrow & \downarrow\\
4 & ... & j+1 & 3 & j+2 & ... & k+3 & 1 & 2 & k+4
\end{matrix}\,,$$

we infer the expression

$$\begin{aligned}
&\rho_{k,j}(4,...,j+1,3,j+2,...,k+3,1,2,k+4)\\
&= \sum_{I\subseteq\{4,...,j+1\}} R^{(|I|)}(I,3,j+2,R^{(k+1-j+|\complement I|)}(\complement I,j+3,...,k+3,1,2,k+4))\\
&- \sum_{I\subseteq\{4,...,j+1\}} \sum_{h=j+3}^{k+3} R^{(k+1-j+|I|)}(I,j+3,...,R^{(|\complement I|)}(\complement I,3,j+2,h),...,k+3,1,2,k+4)\\
&- \sum_{I\subseteq\{4,...,j+1\}} R^{(k+1-j+|I|)}(I,j+3,...,k+3,R^{(|\complement I|)}(\complement I,3,j+2,1),2,k+4)\\
&- \sum_{I\subseteq\{4,...,j+1\}} R^{(k+1-j+|I|)}(I,j+3,...,k+3,1,R^{(|\complement I|)}(\complement I,3,j+2,2),k+4)\\
&- \sum_{I\subseteq\{4,...,j+1\}} R^{(k+1-j+|I|)}(I,j+3,...,k+3,1,2,R^{(|\complement I|)}(\complement I,3,j+2,k+4))\,.
\end{aligned}$$

Performing now the change of variables

$$\begin{matrix}
1 & ... & s-3 & s-2 & s-1 & ... & j-2 & j-1 & j & ... & k+2 & k+3 & k+4\\
\downarrow & & \downarrow & \downarrow & \downarrow & & \downarrow & \downarrow & \downarrow & & \downarrow & \downarrow & \downarrow\\
4 & ... & s & 3 & s+1 & ... & j & 2 & j+1 & ... & k+3 & 1 & k+4
\end{matrix}\,,$$

we infer the expression

$$\begin{aligned}
&\rho_{k,j}(4,...,s,3,s+1,...,j,2,j+1,...,k+3,1,k+4)\\
&= \sum_{I\subseteq\{4,...,s,3,s+1,...,j\}} R^{(|I|)}(I,2,j+1,R^{(k+1-j+|\complement I|)}(\complement I,j+2,...,k+3,1,k+4))\\
&- \sum_{I\subseteq\{4,...,s,3,s+1,...,j\}} \sum_{h=j+2}^{k+3} R^{(k+1-j+|I|)}(I,j+2,...,R^{(|\complement I|)}(\complement I,2,j+1,h),...,k+3,1,k+4)
\end{aligned}$$

$$- \sum_{I\subseteq\{4,...,s,3,s+1,...,j\}} R^{(k+1-j+|I|)}(I, j+2, ..., k+3, R^{(|\complement I|)}(\complement I, 2, j+1, 1), k+4)$$

$$- \sum_{I\subseteq\{4,...,s,3,s+1,...,j\}} R^{(k+1-j+|I|)}(I, j+2, ..., k+3, 1, R^{(|\complement I|)}(\complement I, 2, j+1, k+4))\,.$$

We perform finally the change of variables

$$\begin{array}{cccccccccccc} 1 & ... & j-2 & j-1 & j & ... & s-2 & s-1 & s & ... & k+2 & k+3 & k+4 \\ \downarrow & & \downarrow & \downarrow & \downarrow & & \downarrow & \downarrow & \downarrow & & \downarrow & \downarrow & \downarrow \\ 4 & ... & j+1 & 2 & j+2 & ... & s & 3 & s+1 & ... & k+3 & 1 & k+4 \end{array}\,.$$

We distinguish two cases. In the first case, when $s = j+1$, we have

$$\rho_{k,j}(4, ..., j+1, 2, 3, j+2, ..., k+3, 1, k+4)$$

$$= \sum_{I\subseteq\{4,...,j+1\}} R^{(|I|)}(I, 2, 3, R^{(k+1-j+|\complement I|)}(\complement I, j+2, ..., k+3, 1, k+4))$$

$$- \sum_{\substack{I\subseteq\{4,...,j+1\}\\ j+2\leqslant h\leqslant k+3}} R^{(k+1-j+|I|)}(I, j+2, ..., R^{(|\complement I|)}(\complement I, 2, 3, h), ..., k+3, 1, k+4)$$

$$- \sum_{I\subseteq\{4,...,j+1\}} R^{(k+1-j+|I|)}(I, j+2, ..., k+3, R^{(|\complement I|)}(\complement I, 2, 3, 1), k+4)$$

$$- \sum_{I\subseteq\{4,...,j+1\}} R^{(k+1-j+|I|)}(I, j+2, ..., k+3, 1, R^{(|\complement I|)}(\complement I, 2, 3, k+4))\,,$$

and in the second case, when $j+2 \leqslant s \leqslant k+3$, we have

$$\rho_{k,j}(4, ..., j+1, 2, j+2, ..., s, 3, s+1, ..., k+3, 1, k+4)$$

$$= \sum_{I\subseteq\{4,...,j+1\}} R^{(|I|)}(I, 2, j+2, R^{(k+1-j+|\complement I|)}(\complement I, j+3, ..., s, 3, s+1, ..., k+3, 1, k+4))$$

$$- \sum_{\substack{I\subseteq\{4,...,j+1\}\\ j+3\leqslant h\leqslant s}} R^{(k+1-j+|I|)}(I, j+3, ..., R^{(|\complement I|)}(\complement I, 2, j+2, h), ..., s, 3, s+1, ..., k+3, 1, k+4)$$

$$- \sum_{I\subseteq\{4,...,j+1\}} R^{(k+1-j+|I|)}(I, j+3, ..., s, R^{(|\complement I|)}(\complement I, 2, j+2, 3), s+1, ..., k+3, 1, k+4)$$

$$- \sum_{\substack{I\subseteq\{4,...,j+1\}\\ s+1\leqslant h\leqslant k+3}} R^{(k+1-j+|I|)}(I, j+3, ..., s, 3, s+1, ..., R^{(|\complement I|)}(\complement I, 2, j+2, h), ..., k+3, 1, k+4)$$

$$- \sum_{I\subseteq\{4,...,j+1\}} R^{(k+1-j+|I|)}(I, j+3, ..., s, 3, s+1, ..., k+3, R^{(|\complement I|)}(\complement I, 2, j+2, 1), k+4)$$

$$- \sum_{I\subseteq\{4,...,j+1\}} R^{(k+1-j+|I|)}(I, j+3, ..., s, 3, s+1, ..., k+3, 1, R^{(|\complement I|)}(\complement I, 2, j+2, k+4))\,.$$

Then the conclusion follows from lemma 4.2. □

## 4.2 Expliciting the 1-exterior product of the lower order term of $S_k$

Given a subset $S \subset \mathbb{Z}_{\geqslant 1}$ and $I \subset S$ we denote below by abuse of notation $\complement I \subset S$ the complementary set of $I$ inside $S$. We show now the following key result.

**Proposition 4.4.** *For all $h \geqslant 2$, let*

$$\hat{S}_h \;:=\; \frac{2 i^{h+1}}{(h+1)!\,h!} \mathrm{Sym}_{2,\ldots,h+1} \Phi_{h-3}\,.$$

*Then for all $k \geqslant 1$, hold the equality*

$$\mathrm{Sym}_{3,\ldots,k+4}\left[\sum_{r=3}^{k+2} (r+1)! \sum_{p=2}^{r-1} (i d_1^{\nabla})^{k+2-r} (p \hat{S}_p \wedge_1 \hat{S}_{r-p+1})\right] \;=\; \mathrm{Sym}_{3,\ldots,k+4} \mathbb{L}(k)\,,$$

*with*

$$\mathbb{L}(k) \;:=\; \sum_{r=3}^{k+2} \sum_{p=2}^{r-1} \frac{4\,(-i)^{k+1} (r+1)!\, \mathbb{L}_{k,r,p}}{(r-p+2)!\,(p+1)!}\,,$$

*and with*

$$\begin{aligned}
& \mathbb{L}_{k,r,p}(1,\ldots,k+4) \\
:= \; & \sum_{I \subseteq \{1,3,\ldots,k-r+3\}} \\
& R^{(p+|I|-2)}(I, k_p+2, \ldots, k+3, R^{(r-p+|\complement I|-1)}(\complement I, k_r+1, \ldots, k_p, 2, k_p+1), k_r, k+4) \\
- \; & \sum_{I \subseteq \{2,\ldots,k-r+3\}} \\
& R^{(p+|I|-2)}(I, k_p+2, \ldots, k+3, R^{(r-p+|\complement I|-1)}(\complement I, k_r+1, \ldots, k_p, 1, k_p+1), k_r, k+4)\,.
\end{aligned}$$

*with $k_p := k-p+4$, in the case $k+2-r \geqslant 1$ and with*

$$\begin{aligned}
& \mathbb{L}_{k,k+2,p}(1,\ldots,k+4) \\
:= \; & \sum_{j=k_p+1}^{k+3} R^{(p-2)}(k_p+2, \ldots, j, R^{(k-p+1)}(3, \ldots, k_p, 2, k_p+1), j+1, \ldots, k+3, 1, k+4) \\
+ \; & R^{(p-2)}(k_p+2, \ldots, k+4, 1, R^{(k-p+1)}(3, \ldots, k_p, 2, k_p+1))
\end{aligned}$$

$$- \sum_{j=k_p+1}^{k+3} R^{(p-2)}(k_p+2, ..., j, R^{(k-p+1)}(3, ..., k_p, 1, k_p+1), j+1, ..., k+3, 2, k+4)$$

$$- \; R^{(p-2)}(k_p+2, ..., k+4, 2, R^{(k-p+1)}(3, ..., k_p, 1, k_p+1)) \, .$$

**Proof.** We notice that in the case $k+2-r \geqslant 1$,

$$\mathrm{Sym}_{3,...,k+4} [(i d_1^{\nabla})^{k+2-r}(p \, \hat{S}_p \wedge_1 \hat{S}_{r-p+1})] \; = \; (-i)^{k+2-r} \, \mathrm{Sym}_{3,...,k+4}(\varphi_{k,r,p} - \psi_{k,r,p}) \, ,$$

with

$$\varphi_{k,r,p}(1, ..., k+4) \; := \; \nabla^{k+2-r}(p \, \hat{S}_p \wedge_1 \hat{S}_{r-p+1})(2, ..., k+3-r, 1, k+4-r, ..., k+4) \, ,$$

$$\psi_{k,r,p}(1, ..., k+4) \; := \; \nabla^{k+2-r}(p \, \hat{S}_p \wedge_1 \hat{S}_{r-p+1})(1, 3, ..., k+3-r, 2, k+4-r, ..., k+4) \, .$$

We notice now that

$$\mathrm{Sym}_{3,...,k+4} \, \varphi_{k,r,p} \; = \; \mathrm{Sym}_{3,...,k+4} \, \tilde{\varphi}_{k,r,p} \, ,$$

$$\mathrm{Sym}_{3,...,k+4} \, \psi_{k,r,p} \; = \; \mathrm{Sym}_{3,...,k+4} \, \tilde{\psi}_{k,r,p} \, ,$$

with

$$\tilde{\varphi}_{k,r,p}(1, ..., k+4) \; := \; \nabla^{k+2-r} \, \lambda_{r,p} \, (2, ..., k+3-r, 1, k+4-r, ..., k+4) \, ,$$

$$\tilde{\psi}_{k,r,p}(1, ..., k+4) \; := \; \nabla^{k+2-r} \lambda_{r,p} \, (1, 3, ..., k+3-r, 2, k+4-r, ..., k+4) \, ,$$

with

$$\lambda_{r,p} \, (1, k+4-r, ..., k+4)$$

$$:= \; -p \, \hat{S}_p(k+4-r, \hat{S}_{r-p+1}(1, k+5-r, ..., k+5-p), k+6-p, ..., k+4) \, .$$

Therefore

$$(\tilde{\varphi}_{k,r,p} - \tilde{\psi}_{k,r,p})(1, ..., k+4)$$

$$= \; p \sum_{I \subseteq \{1,3,...,k-r+3\}} \hat{S}_p^{(|I|)}\big(I, k_r, \hat{S}_{r-p+1}^{(|\complement I|)}(\complement I, 2, k_r+1, ..., k_p+1), k_p+2, ..., k+4\big)$$

$$- \; p \sum_{I \subseteq \{2,...,k-r+3\}} \hat{S}_p^{(|I|)}\big(I, k_r, \hat{S}_{r-p+1}^{(|\complement I|)}(\complement I, 1, k_r+1, ..., k_p+1), k_p+2, ..., k+4\big) \, .$$

We notice now that for any tensor $\Phi \equiv \Phi(1, ..., p+1)$, the tensor

$$\tilde{\Phi}(1, ..., p+1) \; := \; \sum_{j=2}^{p+1} \Phi(1, 3, ..., j, 2, j+1, ..., p+1) \, ,$$

satisfies the identity

$$\mathrm{Sym}_{2,...,p+1} \, \Phi \; = \; \mathrm{Sym}_{3,...,p+1} \, \tilde{\Phi} \, .$$

Thus

$$(\mathrm{Sym}_{2,...,p+1}\,\Phi)^{(h)} \;=\; \mathrm{Sym}_{h+3,...,h+p+1}\,\tilde{\Phi}^{(h)}\,.$$

We infer

$$\begin{aligned}
&(\tilde{\varphi}_{k,r,p} - \tilde{\psi}_{k,r,p})(1,...,k+4)\\
=\;& \frac{4i^{r+3}}{(p+1)!(p-1)!(r-p+2)!(r-p+1)!}\sum_{I\subseteq\{1,3,...,k-r+3\}}\\
&\big(\mathrm{Sym}_{|I|+3,...,|I|+p+1}\tilde{\Phi}^{(|I|)}_{p-3}\big)\big(I,\,k_r,\,\big(\mathrm{Sym}_{|\complement I|+2,...,|\complement I|+r-p+2}\Phi^{(|\complement I|)}_{r-p-2}\big)(\complement I,\,2,\,k_r+1,\,...,\,k_p+1),\\
&k_p+2,...,k+4\big)\\
-\;& \frac{4i^{r+3}}{(p+1)!(p-1)!(r-p+2)!(r-p+1)!}\sum_{I\subseteq\{2,...,k-r+3\}}\\
&\big(\mathrm{Sym}_{|I|+3,...,|I|+p+1}\tilde{\Phi}^{(|I|)}_{p-3}\big)\big(I,\,k_r,\,\big(\mathrm{Sym}_{|\complement I|+2,...,|\complement I|+r-p+2}\Phi^{(|\complement I|)}_{r-p-2}\big)(\complement I,\,1,\,k_r+1,\,...,\,k_p+1),\\
&k_p+2,...,k+4\big)\,.
\end{aligned}$$

and thus

$$(-i)^{k+2-r}\,\mathrm{Sym}_{3,...,k+4}(\tilde{\varphi}_{k,r,p} - \tilde{\psi}_{k,r,p}) \;=\; \frac{4\,(-1)^{k-r}\,i^{k+1}}{(r-p+2)!\,(p+1)!}\,\mathrm{Sym}_{3,...,k+4}\,\tilde{\mathbb{L}}_{k,r,p}\,,$$

with

$$\begin{aligned}
&\tilde{\mathbb{L}}_{k,r,p}(1,...,k+4)\\
:=\;& \sum_{I\subseteq\{1,3,...,k-r+3\}}\tilde{\Phi}^{(|I|)}_{p-3}\big(I,k_r,\Phi^{(|\complement I|)}_{r-p-2}(\complement I,2,k_r+1,...,k_p+1),k_p+2,...,k+4\big)\\
-\;& \sum_{I\subseteq\{2,...,k-r+3\}}\tilde{\Phi}^{(|I|)}_{p-3}\big(I,k_r,\Phi^{(|\complement I|)}_{r-p-2}(\complement I,1,k_r+1,...,k_p+1),k_p+2,...,k+4\big)\\
=\;& (-1)^{r-p-1}\sum_{I\subseteq\{1,3,...,k-r+3\}}\\
&\tilde{\Phi}^{(|I|)}_{p-3}(I,k_r,R^{(r-p+|\complement I|-1)}(\complement I,k_r+1,...,k_p,2,k_p+1),k_p+2,...,k+4)\\
-\;& (-1)^{r-p-1}\sum_{I\subseteq\{2,...,k-r+3\}}\\
&\tilde{\Phi}^{(|I|)}_{p-3}(I,k_r,R^{(r-p+|\complement I|-1)}(\complement I,k_r+1,...,k_p,1,k_p+1),k_p+2,...,k+4)\,.
\end{aligned}$$

Using the expression

$$\begin{aligned}
\tilde{\Phi}_{p-3}(1,...,p+1) \;=\;& (-1)^p\sum_{j=2}^{p} R^{(p)}(3,...,j,2,j+1,...,p,1,p+1)\\
&+\;(-1)^pR^{(p)}(3,...,p+1,1,2)\,, \qquad (4.3)
\end{aligned}$$

we deduce the formula

$$(-i)^{k+2-r}\,\mathrm{Sym}_{3,\ldots,k+4}(\tilde{\varphi}_{k,r,p}-\tilde{\psi}_{k,r,p}) \;=\; \frac{4\,(-i)^{k+1}}{(r-p+2)!\,(p+1)!}\,\mathrm{Sym}_{3,\ldots,k+4}\,\hat{\mathbb{L}}_{k,r,p}\,,$$

with

$$\begin{aligned}
&\hat{\mathbb{L}}_{k,r,p}(1,\ldots,k+4)\\
:=\;& \sum_{I\subseteq\{1,3,\ldots,k-r+3\}}\;\sum_{j=k_p+1}^{k+3}\\
& R^{(p+|I|-2)}(I,k_p+2,\ldots,j,R^{(r-p+|\complement I|-1)}(\complement I,k_r+1,\ldots,k_p,2,k_p+1),j+1,\ldots,k+3,k_r,k+4)\\
+\;& \sum_{I\subseteq\{1,3,\ldots,k-r+3\}} R^{(p+|I|-2)}(I,k_p+2,\ldots,k+4,k_r,R^{(r-p+|\complement I|-1)}(\complement I,k_r+1,\ldots,k_p,2,k_p+1))\\
-\;& \sum_{I\subseteq\{2,\ldots,k-r+3\}}\;\sum_{j=k_p+1}^{k+3}\\
& R^{(p+|I|-2)}(I,k_p+2,\ldots,j,R^{(r-p+|\complement I|-1)}(\complement I,k_r+1,\ldots,k_p,1,k_p+1),j+1,\ldots,k+3,k_r,k+4)\\
-\;& \sum_{I\subseteq\{2,\ldots,k-r+3\}} R^{(p+|I|-2)}(I,k_p+2,\ldots,k+4,k_r,R^{(r-p+|\complement I|-1)}(\complement I,k_r+1,\ldots,k_p,1,k_p+1))\,.
\end{aligned}$$

The fact that in the case under consideration $k+2-r\geqslant 1$, i.e $k_r\geqslant 3$, implies

$$\mathrm{Sym}_{3,\ldots,k+4}\left[(i\,d_1^{\nabla})^{k+2-r}(p\,\hat{S}_p\wedge_1\hat{S}_{r-p+1})\right] \;=\; \frac{4\,(-i)^{k+1}}{(r-p+2)!\,(p+1)!}\,\mathrm{Sym}_{3,\ldots,k+4}\,\mathbb{L}_{k,r,p}\,.$$

We consider now the case $r=k+2$. We consider the expansion

$$\begin{aligned}
&(p\,\hat{S}_p\wedge_1\hat{S}_{k-p+3})(1,\ldots,k+4)\\
=\;& p\,\hat{S}_p(1,\hat{S}_{k-p+3}(2,\ldots,k_p+1),k_p+2,\ldots,k+4)\\
-\;& p\,\hat{S}_p(2,\hat{S}_{k-p+3}(1,3,\ldots,k_p+1),k_p+2,\ldots,k+4)\\
=\;& \frac{4i^{k+5}}{(p+1)!(p-1)!(k-p+4)!(k-p+3)!}\times\\
\times\;& (\mathrm{Sym}_{3,\ldots,p+1}\tilde{\Phi}_{p-3})(1,(\mathrm{Sym}_{2,\ldots,k_p}\Phi_{k-p})(2,\ldots,k_p+1),k_p+2,\ldots,k+4)\\
-\;& \frac{4i^{k+5}}{(p+1)!(p-1)!(k-p+4)!(k-p+3)!}\times\\
\times\;& (\mathrm{Sym}_{3,\ldots,p+1}\tilde{\Phi}_{p-3})(2,(\mathrm{Sym}_{2,\ldots,k_p}\Phi_{k-p})(1,3,\ldots,k_p+1),k_p+2,\ldots,k+4)\,.
\end{aligned}$$

We infer

$$\mathrm{Sym}_{3,\ldots,k+4}\,(p\,\hat{S}_p\wedge_1\hat{S}_{r-p+1}) \;=\; \frac{4i^{k+5}}{(p+1)!(k-p+4)!}\,\mathrm{Sym}_{3,\ldots,k+4}\,\tilde{\mathbb{L}}_{k,k+2,p}\,,$$

with

$$
\begin{aligned}
& \tilde{\mathbb{L}}_{k,k+2,p}(1, ..., k+4) \\
= \;& \tilde{\Phi}_{p-3}(1, \Phi_{k-p}(2, ..., k_p+1), k_p+2, ..., k+4) \\
- \;& \tilde{\Phi}_{p-3}(2, \Phi_{k-p}(1, 3, ..., k_p+1), k_p+2, ..., k+4) \\
= \;& (-1)^{k-p+1}\, \tilde{\Phi}_{p-3}(1, R^{(k-p+1)}(3, ..., k_p, 2, k_p+1), k_p+2, ..., k+4) \\
- \;& (-1)^{k-p+1}\, \tilde{\Phi}_{p-3}(2, R^{(k-p+1)}(3, ..., k_p, 1, k_p+1), k_p+2, ..., k+4)\,.
\end{aligned}
$$

Using (4.4), we obtain

$$
\mathrm{Sym}_{3,...,k+4}\,(p\, \hat{S}_p \wedge_1 \hat{S}_{r-p+1}) \;=\; \frac{4\,(-i)^{k+1}}{(p+1)!(k-p+4)!}\,\mathrm{Sym}_{3,...,k+4}\,\mathbb{L}_{k,k+2,p}\,,
$$

which shows the required formula. □

# 5 Vanishing of the integrability conditions in the cases $k = 4, 5$

In this section we want to show the identity

$$
\mathrm{Circ}\,\mathrm{Sym}_{3,...,k+4}\left[2i(id)^k(R^{'})_2 - \sum_{r=3}^{k+2}(r+1)!\sum_{p=2}^{r-1}(id)^{k+2-r}(p\,\hat{S}_p \wedge \hat{S}_{r-p+1})\right] \;=\; 0,
$$

for all integers $k = 1, 2$. We notice that $\hat{S}_p = S_p$ for $p = 2, 3$. This corresponds to the cases $k = 4, 5$ in the main integrability conditions $\mathrm{Circ}\,\mathrm{Sym}_{3,...,k+1}\,\Theta_k^{\nabla} = 0$.

## 5.1 Alternative proof of the vanishing of the integrability conditions in the case $k = 4$

We set for simplicity $ab, cde := R(a, b, R(c, d, e))$ and a similar definition for $abc, de$ or $a, bcd, e$.

We use here the notation $R^{'} := \nabla R^{\nabla}$. We observe first that lemma 4.3 implies the identity

$$
\begin{aligned}
\mathrm{Circ}\,\mathrm{Sym}_{3,4,5}[d\,(R^{'})_2] \;&=\; \mathrm{Circ}\,\mathrm{Sym}_{4,5}\,\mathrm{Alt}_2\,\mathbb{Q}_{1,2} \\
\mathbb{Q}_{1,2}(1, ..., 5) \;&=\; 23, 415 - 234, 15 - 41, 235 \\
&+\; 24, 315 - 243, 15 - 3, 241, 5 - 31, 245\,.
\end{aligned}
$$

We notice that $\mathrm{Circ}\,\mathbb{Q}_{1,2} = \mathrm{Circ}\,\tilde{\mathbb{Q}}_{1,2}$, with

$$
\begin{aligned}
\tilde{\mathbb{Q}}_{1,2}(1, ..., 5) \;&=\; 12, 435_1 - 124, 35_2 - 43, 125_3 \\
&+\; 34, 125_3 - 142, 35_2 - 3, 241, 5_2 - 12, 345_1 \\
&=\; -2\cdot 12, 345 + 2\cdot 34, 125 - 124, 35 - 142, 35 - 421, 35 \\
&=\; -2\cdot 12, 345 + 2\cdot 34, 125 - 124, 35 + 214, 35 \\
&=\; -2\cdot 12, 345 + 2\cdot 34, 125 - 2\cdot 124, 35\,.
\end{aligned}
$$

(We denote by a subscript $j=1,2,3$, the terms that we combine together). Using the fact that all the operators $\mathrm{Circ}, \mathrm{Sym}_{4,5}, \mathrm{Alt}_2$ commute we infer the formula

$$\mathrm{Circ}\,\mathrm{Sym}_{3,4,5}[d\,(R')_2] \;=\; -4\,\mathrm{Circ}\,\mathrm{Sym}_{4,5}\,\boldsymbol{Q}_1$$

$$\boldsymbol{Q}_1(1,...,5) \;=\; 12,345-3,124,5-34,125\,.$$

On the other hand using proposition 4.4, we infer the formula

$$S_2 \;:=\; -\frac{i^3}{6}\mathrm{Sym}_{2,3}\,R\,,$$

$$4!\,2\,\mathrm{Sym}_{3,4,5}(S_2\wedge S_2) \;=\; -\frac{8}{3}\mathrm{Sym}_{3,4,5}\,\mathrm{Alt}_2\,\mathbb{L}_{1,3,2}\,,$$

$$\mathbb{L}_{1,3,2}(1,...,5) \;:=\; 324,15+15,234\,.$$

We show now the identity

$$\mathrm{Circ}\,\mathrm{Sym}_{3,4,5}\big[\,2\,d\,(R')_2+4!\,2\,\mathrm{Sym}_{3,4,5}(S_2\wedge S_2)\big] \;=\; 0\,,$$

i.e.

$$12\,\mathrm{Circ}\,\mathrm{Sym}_{4,5}\,\boldsymbol{Q}_1-4\,\mathrm{Circ}\,\mathrm{Sym}_{3,4,5}\mathrm{Alt}_2\,\mathbb{L}_{1,3,2}=0\,. \tag{5.1}$$

Using the identity $[\mathrm{Sym}_{3,4,5}, \mathrm{Alt}_2]=0$, the equalities

$$\begin{aligned}
\mathrm{Sym}_{3,4,5}\,\mathbb{L}_{1,3,2} &= \mathrm{Sym}_{4,5}\,[324,15+423,15+425,13+14,235+13,245+14,253]\\
&= \mathrm{Sym}_{4,5}\,[324,15+423,15+425,13+14,235+31,425+14,253]\,,\\
&\quad \mathrm{Circ}\,\mathrm{Sym}_{4,5}[324,15+423,15+425,13+14,235+23,415+14,253]\\
&= \mathrm{Circ}\,\mathrm{Sym}_{4,5}[214,35+412,35+435,21+34,125+12,435+34,152]\\
&= \mathrm{Circ}\,\mathrm{Sym}_{4,5}[214,35+412,35-345,21+34,125-12,345+34,152]\,,\\
&\quad \mathrm{Alt}_2[214,35+412,35-345,21+34,125-12,345+34,152]\\
&= 214,35+412,35-345,21+34,125-12,345+34,152\\
&- \;124,35-421,35+345,12-34,215+21,345-34,251\\
&= -2\cdot 124,35+214,35-12,345+2\cdot 34,125+2\cdot 21,345+34,125\\
&= -3\cdot 124,35-3\cdot 12,345+3\cdot 34,125\\
&= 3\cdot 3,124,5-3\cdot 12,345+3\cdot 34,125\,.
\end{aligned}$$

We infer

$$\operatorname{Circ}\operatorname{Sym}_{3,4,5}\operatorname{Alt}_2\mathbb{L}_{1,3,2} \;=\; -3\cdot\operatorname{Circ}\operatorname{Sym}_{4,5}\boldsymbol{Q}_1\,,$$

and thus (5.1).

## 5.2 Vanishing of the integrability conditions in the case $k=5$

We want to show the identity

$$\operatorname{Circ}\operatorname{Sym}_{3,\ldots,6}\left[2id^2(R')_2+\sum_{r=3}^{4}(r+1)!\sum_{p=2}^{r-1}(id)^{4-r}(pS_p\wedge S_{r-p+1})\right] \;=\; 0\,.$$

We set for simplicity $ab(cdef) := R(a,\,b,\,R'(c,\,d,\,e,\,f))$ and a similar definition for $abc(def)$. Using lemma 4.3 we infer the identities

$$\operatorname{Circ}\operatorname{Sym}_{3,\ldots,6}\Big[d^2\,(R')_2\Big] \;=\; \operatorname{Circ}\operatorname{Sym}_{4,5,6}\big[\mathbb{T}_{2,3}-\operatorname{Alt}_2(\mathbb{Q}_{2,2}+\mathbb{Q}_{2,3})\big]\,,$$

$$\begin{aligned}
\mathbb{T}_{2,3}(1,\ldots,6) \;=\;& 35(4126)_1+435(126)_2\\
-\;& (4351)26_3-4(351)26_4\\
+\;& (4352)16_3+4(352)16_4\\
-\;& 12(4356)_5-412(356)_6\,,\\
\mathbb{Q}_{2,2}(1,\ldots,6) \;=\;& 23(4516)_5-(234)516_7-4(235)16_4-451(236)_2\\
+\;& 24(3516)_8+24(5316)_1-(245)316_9-(243)516_7\\
-\;& 5(243)16_4-3(245)16_{10}-35(241)6_{11}-53(241)6_4\\
-\;& 351(246)_6-531(246)_6\,,\\
\mathbb{Q}_{2,3}(1,\ldots,6) \;=\;& 24(3516)_8+324(516)_6-(3245)16_3-3(245)16_{10}\\
-\;& 5(3241)6_{12}-35(241)6_{11}-51(3246)_8-351(246)_6\\
+\;& 23(4516)_5+423(516)_6-(4235)16_3-4(235)16_4\\
-\;& 51(4236)_8-451(236)_2+25(4316)_1+425(316)_2\\
-\;& (4253)16_3-4(253)16_4-3(4251)6_3-43(251)6_4\\
-\;& 31(4256)_5-431(256)_6\,,
\end{aligned}$$

where we denote by a subscript $j = 1, \ldots, 12$, the terms that we combine together inside the expression $\operatorname{Circ}\operatorname{Sym}_{4,\ldots,6}[\mathbb{T}_{2,3} - \operatorname{Alt}_2(\mathbb{Q}_{2,2} + \mathbb{Q}_{2,3})]$. ($j = 12$ is combined with itself). Indeed combining:

for $j = 1$,

$$\begin{aligned} & \operatorname{Sym}_{4,5,6}\big[35(4126) - 24(5316) + 14(5326) - 25(4316) + 15(4326)\big] \\ = \ & \operatorname{Sym}_{4,5,6}\big[34(5126) - 2\cdot 24(5316) + 2\cdot 14(5326)\big], \end{aligned}$$

and thus

$$\begin{aligned} & \operatorname{Circ}\operatorname{Sym}_{4,5,6}\big[35(4126) - 24(5316) + 14(5326) - 25(4316) + 15(4326)\big] \\ = \ & \operatorname{Circ}\operatorname{Sym}_{4,5,6}\big[34(5126) - 2\cdot 34(5126) + 2\cdot 34(5216)\big] \\ = \ & -3\operatorname{Circ}\operatorname{Sym}_{4,5,6}[34(5126)], \end{aligned}$$

for $j = 2$,

$$\begin{aligned} & 435(126) + 2\cdot 451(236) - 2\cdot 452(136) - 425(316) + 415(326) \\ = \ & 435(126) + 2\cdot 451(236) - 2\cdot 452(136) - 452(136) + 451(236) \\ = \ & -453(126) + 3\cdot 451(236) - 3\cdot 452(136), \end{aligned}$$

and thus

$$\begin{aligned} & \operatorname{Circ}\big[435(126) + 2\cdot 451(236) - 2\cdot 452(136) - 425(316) + 415(326)\big] \\ = \ & \operatorname{Circ}\big[-453(126) + 3\cdot 453(126) - 3\cdot 453(216)\big] \\ = \ & 5\operatorname{Circ}[453(126)], \end{aligned}$$

for $j = 3$,

$$\begin{aligned} & \operatorname{Circ}\big[(4352)16 - (4351)26 + (3245)16 - (3145)26 \\ + \ & (4235)16 - (4135)26 + (4253)16 - (4153)26 + 3(4251)6 - 3(4152)6\big] \\ = \ & \operatorname{Circ}\big[(4352)16 - (4351)26 + (2145)36 - (1245)36 \\ + \ & (4235)16 - (4325)16 + (4253)16 + (4513)26 - (4251)36 + (4152)36\big] \\ = \ & \operatorname{Circ}\big[(4352)16 - (4351)26 - (2415)36 - (1245)36 \\ + \ & (4235)16 + (4235)16 + (4253)16 + (4532)16 - (4352)16 + (4351)26\big] \end{aligned}$$

$$
\begin{aligned}
&= \operatorname{Circ}\big[\,2\cdot(4235)16+(4125)36+(4253)16+(4532)16\,\big] \\
&= \operatorname{Circ}\big[\,2\cdot(4235)16+(4125)36-(4325)16\,\big] \\
&= \operatorname{Circ}[3(4235)16+(4125)36] \\
&= 4\operatorname{Circ}[(4125)36]\,,
\end{aligned}
$$

for $j=4$,

$$
\begin{aligned}
&\operatorname{Circ}\big[\,4(352)16-4(351)26+2\cdot4(235)16-2\cdot4(135)26 \\
+&\ 4(253)16-4(153)26+43(251)6-43(152)6\,\big] \\
=&\ \operatorname{Circ}\big[\,4(251)36-4(152)36+2\cdot4(125)36-2\cdot4(215)36 \\
+&\ 4(152)36-4(251)36-4(251)36+4(152)36\,\big] \\
=&\ \operatorname{Circ}\big[\,4\cdot4(125)36+4(521)36+4(152)36\,\big] \\
=&\ \operatorname{Circ}\big[\,4\cdot4(125)36-4(215)36\,\big] \\
=&\ 5\operatorname{Circ}[4(125)36]\,.
\end{aligned}
$$

Moreover for the remaining terms we use the equalities

$$
\begin{aligned}
&\operatorname{Sym}_{4,5,6}\big[\,5(243)16-5(143)26+53(241)6-53(142)6\,\big] \\
=&\ \operatorname{Sym}_{4,5,6}\big[\,5(243)16-5(143)26-5(241)36+5(142)36\,\big] \\
=&\ \operatorname{Sym}_{4,5,6}\big[\,4(253)16-4(153)26-4(251)36+4(152)36\,\big]\,,
\end{aligned}
$$

and

$$
\begin{aligned}
&\operatorname{Circ}\operatorname{Sym}_{4,5,6}\big[\,5(243)16-5(143)26+53(241)6-53(142)6\,\big] \\
=&\ \operatorname{Circ}\operatorname{Sym}_{4,5,6}\big[\,4(152)36-4(251)36-4(251)36+4(152)36\,\big] \\
=&\ 2\operatorname{Circ}\operatorname{Sym}_{4,5,6}\big[\,4(152)36+4(521)36\,\big] \\
-&\ 2\operatorname{Circ}\operatorname{Sym}_{4,5,6}[4(215)36] \\
=&\ 2\operatorname{Circ}\operatorname{Sym}_{4,5,6}[4(125)36]\,.
\end{aligned}
$$

Thus for $j=4$ we have

$$\begin{aligned} & \mathrm{Circ}\,\mathrm{Sym}_{4,5,6}\big[4(352)16-4(351)26+2\cdot 4(235)16-2\cdot 4(135)26+4(253)16-4(153)26 \\ + & \ 43(251)6-43(152)6+5(243)16-5(143)26+53(241)6-53(142)6\big] \\ = & \ 7\,\mathrm{Circ}\,\mathrm{Sym}_{4,5,6}[4(125)36]\,. \end{aligned}$$

For $j=5$,

$$\begin{aligned} & \mathrm{Circ}\,\big[-12(4356)-2\cdot 23(4516)+2\cdot 13(4526)+31(4256)-32(4156)\big] \\ = & \ \mathrm{Circ}\,\big[-12(4356)-2\cdot 12(4536)+2\cdot 21(4536)+12(4356)-21(4356)\big] \\ = & \ \mathrm{Circ}\big[-4\cdot 12(4536)+12(4356)\big] \\ = & \ 5\,\mathrm{Circ}[12(4356)]\,. \end{aligned}$$

For $j=6$,

$$\begin{aligned} & \mathrm{Circ}\,\big[-412(356)+351(246)+531(246)-324(516)-423(516)+431(256)+351(246)\big] \\ = & \ \mathrm{Circ}\,\big[-412(356)+152(346)+512(346)-214(536)-412(536)+412(356)+152(346)\big]\,, \end{aligned}$$

and

$$\begin{aligned} & \mathrm{Circ}\ \mathrm{Sym}_{4,5,6}\ \big[-412(356)+351(246)+531(246)-324(516)-423(516)+431(256)+351(246)\big] \\ = & \ \mathrm{Circ}\ \mathrm{Sym}_{4,5,6}\ \big[-412(356)+142(356)+412(356)-214(536)-412(536)+412(356)+142(356)\big] \\ = & \ \mathrm{Circ}\,\mathrm{Sym}_{4,5,6}\,\big[-412(356)+2\cdot 142(356)+3\cdot 412(356)+214(356)\big] \\ = & \ \mathrm{Circ}\,\mathrm{Sym}_{4,5,6}\,\big[-412(356)+142(356)-421(356)+3\cdot 412(356)\big] \\ = & \ \mathrm{Circ}\,\mathrm{Sym}_{4,5,6}\,\big[-412(356)+142(356)+4\cdot 412(356)\big]\,. \end{aligned}$$

Moreover

$$\begin{aligned} & \mathrm{Circ}\,\mathrm{Sym}_{4,5,6}\,\big[-412(356)+142(356)-241(356)+4\cdot 412(356)-4\cdot 421(356)\big] \\ = & \ \mathrm{Circ}\,\mathrm{Sym}_{4,5,6}\big[7\cdot 412(356)+142(356)+214(356)\big] \\ = & \ \mathrm{Circ}\,\mathrm{Sym}_{4,5,6}\big[7\cdot 412(356)-421(356)\big] \\ = & \ 8\,\mathrm{Circ}\,\mathrm{Sym}_{4,5,6}[412(356)]\,. \end{aligned}$$

For $j=7$,

$$
\begin{aligned}
& \mathrm{Circ}\big[(234)516-(134)526+(243)516-(143)526\big] \\
= \ & \mathrm{Circ}\big[(124)536-(214)536+(142)536-(241)536\big] \\
= \ & \mathrm{Circ}\big[2\cdot(124)536-(412)536-(241)536\big] \\
= \ & 3\,\mathrm{Circ}[(124)536]\,.
\end{aligned}
$$

For $j=8$,

$$
\begin{aligned}
& \mathrm{Sym}_{4,5,6}\big[-2\cdot 24(3516)+2\cdot 14(3526)+51(3246)-52(3146)+51(4236)-52(4136)\big] \\
= \ & \mathrm{Sym}_{4,5,6}\big[-2\cdot 24(3516)+2\cdot 14(3526)+41(3256)-42(3156)+41(5236)-42(5136)\big] \\
= \ & \mathrm{Sym}_{4,5,6}\big[3\cdot 24(3156)-3\cdot 14(3256)-14(5236)+24(5136)\big]\,,
\end{aligned}
$$

and

$$
\begin{aligned}
& \mathrm{Circ}\,\mathrm{Sym}_{4,5,6}\big[-2\cdot 24(3516)+2\cdot 14(3526)+51(3246)-52(3146)+51(4236)-52(4136)\big] \\
= \ & \mathrm{Circ}\,\mathrm{Sym}_{4,5,6}\big[3\cdot 34(1256)-3\cdot 34(2156)-34(5126)+34(5216)\big] \\
= \ & \mathrm{Circ}\,\mathrm{Sym}_{4,5,6}\big[3\cdot 34(1256)+3\cdot 34(2516)-2\cdot 34(5126)\big] \\
= \ & -5\,\mathrm{Circ}\,\mathrm{Sym}_{4,5,6}[34(5126)]\,.
\end{aligned}
$$

For $j=9$,

$$
\mathrm{Circ}[(245)316] \ = \ \mathrm{Circ}[(345)126]\,,
$$

and thus

$$
\mathrm{Circ}[(245)316-(145)326] \ = \ 2\,\mathrm{Circ}[(345)126]\,.
$$

For $j=10$,

$$
\begin{aligned}
2\,\mathrm{Circ}[3(245)16-3(145)26] \ &= \ 2\,\mathrm{Circ}[1(345)26-2(345)16] \\
&= \ 2\,\mathrm{Circ}[1(345)26+21(345)6] \\
&= \ -2\,\mathrm{Circ}[(345)216] \\
&= \ 2\,\mathrm{Circ}[(345)126]\,.
\end{aligned}
$$

For $j=11$,

$$\begin{aligned} 2\cdot 35(241)6 - 2\cdot 35(142)6 &= 2\cdot 35(241)6 + 2\cdot 35(412)6 \\ &= -2\cdot 35(124)6 . \end{aligned}$$

For $j=12$,

$$\begin{aligned} 5(3241)6 - 5(3142)6 &= 5(3241)6 + 5(3412)6 \\ &= -5(3124)6 . \end{aligned}$$

We infer the formula

$$\begin{aligned} \mathrm{Circ\,Sym}_{3,\ldots,6}\Big[d^2\,(R')_2\Big] &= \mathrm{Circ\,Sym}_{4,5,6}\,\mathbb{T}_2\,, \\ \mathbb{T}_2(1,\ldots,6) &= -3{\cdot}34(5126)_{_1} + 5\cdot 453(126) + 4\cdot(4125)36 \\ &+\; 7\cdot 4(125)36 + 5\cdot 12(4356) + 8\cdot 412(356) \\ &+\; 3\cdot(124)536 - 5\cdot 34(5126)_{_1} + 4\cdot(345)126 \\ &-\; 2\cdot 35(124)6 - 5(3124)6 . \end{aligned}$$

Moreover

$$\begin{aligned} \mathrm{Circ\,Sym}_{4,5,6}\,\mathbb{T}_2 &= \mathrm{Circ\,Sym}_{4,5,6}\,\tilde{\mathbb{T}}_2\,, \\ \tilde{\mathbb{T}}_2(1,\ldots,6) &= -8\cdot 34(5126) + 5\cdot 453(126) + 4\cdot(4125)36 \\ &+\; 7\cdot 5(124)36_{_1} + 5\cdot 12(4356) + 8\cdot 412(356) \\ &+\; 3\cdot(124)536_{_1} + 4\cdot(345)126 - 2\cdot 35(124)6_{_1} . \end{aligned}$$

For $j=1$, we have

$$\begin{aligned} 7\cdot 5(124)36 - 3\cdot(124)356 - 2\cdot 35(124)6 &= 7\cdot 5(124)36 - (124)356 + 2\cdot 5(124)36 \\ &= 9\cdot 5(124)36 - (124)356 . \end{aligned}$$

We conclude

$$\begin{aligned} \mathrm{Circ\,Sym}_{3,\ldots,6}\Big[d^2\,(R')_2\Big] &= \mathrm{Circ\,Sym}_{4,5,6}\,\boldsymbol{T}_2\,, \\ \boldsymbol{T}_2(1,\ldots,6) &= -8\cdot 34(5126) + 5\cdot 453(126) + 4\cdot(4125)36 \\ &+\; 9\cdot 4(125)36 + 5\cdot 12(4356) + 8\cdot 412(356) \\ &-\; (124)356 + 4\cdot(345)126 . \end{aligned}$$

On the other hand using proposition 4.4, we infer the identity

$$\begin{aligned}
&\mathrm{Sym}_{3,\ldots,6}\left[\sum_{r=3}^{4}(r+1)!\sum_{p=2}^{r-1}(i\,d)^{4-r}(pS_p\wedge S_{r-p+1})\right] \\
&= \frac{2i}{3}\mathrm{Sym}_{3,\ldots,6}\,\mathrm{Alt}_2(4\mathbb{L}_{232}+5\mathbb{L}_{242}+5\,\mathbb{L}_{243})\,,
\end{aligned}$$

$$\begin{aligned}
\mathbb{L}_{232}(1,\ldots,6) &= (1425)36+1(425)36\,,\\
\mathbb{L}_{242}(1,\ldots,6) &= (3425)16+61(3425)\,,\\
\mathbb{L}_{243}(1,\ldots,6) &= (324)516+5(324)16+561(324)\,.
\end{aligned}$$

We write

$$\mathrm{Sym}_{3,\ldots,6}\,\mathrm{Alt}_2(4\mathbb{L}_{232}+5\mathbb{L}_{242}+5\,\mathbb{L}_{243}) \;=\; \mathrm{Sym}_{4,5,6}\,\mathrm{Alt}_2\,\mathbb{L}_2\,,$$

$$\begin{aligned}
\mathbb{L}_2(1,\ldots,6) &= 4\cdot(1425)36+4\cdot(1325)46+4\cdot(1423)56+4\cdot(1425)63\\
&+\; 4\cdot1(425)36+4\cdot1(325)46+4\cdot1(423)56+4\cdot1(425)63\\
&+\; 5\cdot(3425)16+5\cdot(4325)16+5\cdot(5423)16+5\cdot(6425)13\\
&+\; 5\cdot61(3425)+5\cdot61(4325)+5\cdot61(5423)+5\cdot31(6425)\\
&+\; 5\cdot(324)516+5\cdot(423)516+5\cdot(524)316+5\cdot(624)513\\
&+\; 5\cdot5(324)16+5\cdot5(423)16+5\cdot3(524)16+5\cdot5(624)13\\
&+\; 5\cdot561(324)+5\cdot561(423)+5\cdot361(524)+5\cdot531(624)\,.
\end{aligned}$$

We notice that $\mathrm{Circ}\,\mathbb{L}_2=\mathrm{Circ}\,\tilde{\mathbb{L}}_2$, with

$$\begin{aligned}
\tilde{\mathbb{L}}_2(1,\ldots,6) &= 4\cdot(1425)36+4\cdot(1423)56+4\cdot(1425)63\\
&+\; 4\cdot2(435)16+4\cdot3(215)46+4\cdot3(412)56+4\cdot2(435)61\\
&+\; 5\cdot(2415)36+5\cdot(4215)36+5\cdot(5412)36+5\cdot(6435)21\\
&+\; 5\cdot63(2415)+5\cdot63(4215)+5\cdot63(5412)+5\cdot12(6435)\\
&+\; 5\cdot(214)536+5\cdot(412)536+5\cdot(534)126+5\cdot(634)521\\
&+\; 5\cdot5(214)36+5\cdot5(412)36+5\cdot1(534)26+5\cdot5(634)21\\
&+\; 5\cdot563(214)+5\cdot563(412)+5\cdot162(534)+5\cdot512(634)\,,
\end{aligned}$$

and $\mathrm{Sym}_{4,5,6}\,\tilde{\mathbb{L}}_2 = \mathrm{Sym}_{4,5,6}\,\boldsymbol{L}_2$, with

$$\begin{aligned}
\boldsymbol{L}_2(1,...,6) \;=\;& 4\cdot(1425)36+4\cdot(1423)56+4\cdot(1425)63\\
+\;& 4\cdot 2(435)16+4\cdot 3(214)56+4\cdot 3(412)56+4\cdot 2(435)61\\
+\;& 5\cdot(2415)36+5\cdot(4215)36+5\cdot(4512)36+5\cdot(6435)21\\
+\;& 5\cdot 63(2415)+5\cdot 63(4215)+5\cdot 63(4512)+5\cdot 12(6435)\\
+\;& 5\cdot(214)536+5\cdot(412)536+5\cdot(435)126+5\cdot(435)621\\
+\;& 5\cdot 5(214)36+5\cdot 5(412)36+5\cdot 1(435)26+5\cdot 5(435)21\\
+\;& 5\cdot 563(214)+5\cdot 563(412)+5\cdot 162(435)+5\cdot 612(435)\,.
\end{aligned}$$

We write now

$$\begin{aligned}
(\mathrm{Alt}_2\,\boldsymbol{L}_2)(1,...,6) \;=\;& 4\cdot(1425)36_{_1}+4\cdot(1423)56_{_3}+4\cdot(1425)63_{_4}\\
-\;& 4\cdot(2415)36_{_1}-4\cdot(2413)56_{_3}-4\cdot(2415)63_{_4}\\
+\;& 4\cdot 2(435)16_{_5}+4\cdot 3(214)56_{_6}+4\cdot 3(412)56_{_6}+4\cdot 2(435)61_{_7}\\
-\;& 4\cdot 1(435)26_{_5}-4\cdot 3(124)56_{_6}-4\cdot 3(421)56_{_6}-4\cdot 1(435)62_{_7}\\
+\;& 5\cdot(2415)36_{_1}+5\cdot(4215)36_{_8}+5\cdot(4512)36_{_8}+5\cdot(6435)21_{_9}\\
-\;& 5\cdot(1425)36_{_1}-5\cdot(4125)36_{_8}-5\cdot(4521)36_{_8}-5\cdot(6435)12_{_9}\\
+\;& 5\cdot 63(2415)_{_{10}}+5\cdot 63(4215)_{_{10}}+5\cdot 63(4512)_{_{10}}+5\cdot 12(6435)_{_9}\\
-\;& 5\cdot 63(1425)_{_{10}}-5\cdot 63(4125)_{_{10}}-5\cdot 63(4521)_{_{10}}-5\cdot 21(6435)_{_9}\\
+\;& 5\cdot(214)536_{_{11}}+5\cdot(412)536_{_{11}}+5\cdot(435)126_{_5}+5\cdot(435)621_{_7}\\
-\;& 5\cdot(124)536_{_{11}}-5\cdot(421)536_{_{11}}-5\cdot(435)216_{_5}-5\cdot(435)612_{_7}\\
+\;& 5\cdot 5(214)36_{_{12}}+5\cdot 5(412)36_{_{12}}+5\cdot 1(435)26_{_5}+5\cdot 6(435)21_{_7}\\
-\;& 5\cdot 5(124)36_{_{12}}-5\cdot 5(421)36_{_{12}}-5\cdot 2(435)16_{_5}-5\cdot 6(435)12_{_7}\\
+\;& 5\cdot 563(214)_{_{13}}+5\cdot 563(412)_{_{13}}+5\cdot 162(435)_{_7}+5\cdot 612(435)_{_7}\\
-\;& 5\cdot 563(124)_{_{13}}-5\cdot 563(421)_{_{13}}-5\cdot 261(435)_{_7}-5\cdot 621(435)_{_7}\,.
\end{aligned}$$

We need to explain the details for the case $j=7$. In this case we sum

$$
\begin{aligned}
& 4\cdot 2(435)61+5\cdot(435)621+5\cdot 6(435)21_{0}+5\cdot 162(435)_{a}+5\cdot 612(435)_{a}\\
- & 4\cdot 1(435)62-5\cdot(435)612-5\cdot 6(435)12_{0}-5\cdot 261(435)_{a}-5\cdot 621(435)_{a}\\
= & 4\cdot 2(435)61_{1}+4\cdot(435)621_{1}+(435)621_{3}+5\cdot 612(435)_{0}+15\cdot 612(435)_{a}\\
- & 4\cdot 1(435)62_{2}-4\cdot(435)612_{2}-(435)612_{3}\\
= & -4\cdot 62(435)1_{1}+4\cdot 61(435)2_{2}+(435)126_{3}+20\cdot 612(435)\\
= & 24\cdot 612(435)+(435)126.
\end{aligned}
$$

The result is (we keep the indices according to the sums on $j=1,\ldots,13$.)

$$
\begin{aligned}
(\mathrm{Alt}_2\,\boldsymbol{L}_2)(1,\ldots,6) = & -(4125)36_{1}+4\cdot(4123)56_{3}+4\cdot(4125)63_{4}\\
+ & 11\cdot(435)126_{5}-12\cdot 3(124)56_{6}+24\cdot 612(435)_{7}+(435)126_{7}\\
- & 15\cdot(4125)36_{8}+15\cdot 12(6435)_{9}-20\cdot 63(4125)_{10}-15\cdot(124)536_{11}\\
- & 15\cdot 5(124)36_{12}-15\cdot 563(124)_{13}.
\end{aligned}
$$

We recombine now the terms of the previous sum under the form

$$
\begin{aligned}
(\mathrm{Alt}_2\,\boldsymbol{L}_2)(1,\ldots,6) = & -16\cdot(4125)36_{1}+4\cdot(4123)56+4\cdot(4125)63_{1}\\
+ & 12\cdot(435)126-12\cdot 3(124)56_{2}+24\cdot 612(435)\\
+ & 15\cdot 12(6435)-20\cdot 63(4125)_{1}-15\cdot(124)536_{2}\\
- & 15\cdot 5(124)36_{2}-15\cdot 563(124).
\end{aligned}
$$

Indeed for $j=1$,

$$
\begin{aligned}
& -16\cdot(4125)36+4\cdot(4125)63-20\cdot 63(4125)\\
= & -12\cdot(4125)36-4\cdot(4125)36+4\cdot(4125)63-20\cdot 63(4125)\\
= & -12\cdot(4125)36+4\cdot 3(4125)6+4\cdot(4125)63-20\cdot 63(4125)\\
= & -12\cdot(4125)36-24\cdot 63(4125),
\end{aligned}
$$

and for $j=2$,

$$
\begin{aligned}
& -12\cdot 3(124)56-15\cdot(124)536-15\cdot 5(124)36\\
= & 12\cdot 53(124)5-3\cdot(124)536-15\cdot 5(124)36\\
= & -27\cdot 5(124)36-3\cdot(124)536.
\end{aligned}
$$

In conclusion

$$\operatorname{Circ}\operatorname{Sym}_{3,\dots,6}\left[\sum_{r=3}^{4}(r+1)!\sum_{p=2}^{r-1}(i d)^{4-r}(p S_p\wedge S_{r-p+1})\right] = \frac{2i}{3}\operatorname{Circ}\operatorname{Sym}_{4,5,6}\tilde{\mathbf{\Lambda}}_2\,,$$

with

$$\begin{aligned}\tilde{\mathbf{\Lambda}}_2(1,\dots,6) &= -12\cdot(4125)36+24\cdot 36(4125)+12\cdot(435)126+24\cdot 612(435)\\ &+ 15\cdot 12(6435)-27\cdot 5(124)36-3\cdot(124)536-15\cdot 563(124)\,.\end{aligned}$$

We notice that $\operatorname{Sym}_{4,5,6}\tilde{\mathbf{\Lambda}}_2=\operatorname{Sym}_{4,5,6}\mathbf{\Lambda}_2$, with

$$\begin{aligned}\mathbf{\Lambda}_2(1,\dots,6) &= -12\cdot(4125)36+24\cdot 34(5126)+12\cdot(435)126+24\cdot 412(536)\\ &+ 15\cdot 12(4536)-27\cdot 4(125)36-3\cdot(124)536-15\cdot 453(126)\,.\end{aligned}$$

We deduce $\mathbf{\Lambda}_3=-3\boldsymbol{T}_2$, which implies the required conclusion in the case $k=2$.

# 6 Appendix

We remind first the following very elementary lemmas. The arguments below can be found in [Pali].

**Lemma 6.1.** *Let $M$ be a smooth manifold, let also $\nabla$ be a covariant derivative operator acting on the smooth sections of $T_M$ and consider the vector field over $T_M$*

$$\zeta_\eta^\nabla := H_\eta^\nabla\cdot\eta\,, \tag{6.1}$$

*$\eta\in T_M$. We denote by $\Phi_t^\nabla$ the corresponding 1-parameter subgroup of transformations of $T_M$. Then for any $\eta\in T_M$ the curve $c_t:=\pi_{T_M}\circ\Phi_t^\nabla(\eta)$ is the $\nabla$-geodesic with initial speed $\dot{c}_0=\eta$ and $\dot{c}_t=\Phi_t^\nabla(\eta)$.*

**Proof.** The flow line $\eta_t:=\Phi_t^\nabla(\eta)$ satisfies the identity

$$\dot{\eta}_t=H_{\eta_t}^\nabla\cdot\eta_t\,. \tag{6.2}$$

We deduce

$$\begin{aligned}\dot{c}_t &= d_{\eta_t}\pi_{T_M}\cdot\dot{\eta}_t\\ &= d_{\eta_t}\pi_{T_M}\cdot H_{\eta_t}^\nabla\cdot\eta_t\\ &= \eta_t\,,\end{aligned}$$

and $\ddot{c}_t=H_{\dot{c}_t}^\nabla\cdot\dot{c}_t$, which is the $\nabla$-geodesic equation. □

**Corollary 6.2.** *Let $M$ be a smooth manifold, let also $\nabla$ be a covariant derivative operator acting on the smooth sections of $T_M$ and let $U$ be an open neighborhood $U$ of $M$ inside $T_M$. A complex structure $J$ over $U$ satisfies the conditions*

$$J_{|M}=J^{\mathrm{can}}\,, \tag{6.3}$$

$$J_\eta H_\eta^\nabla\cdot\eta=T_\eta\cdot\eta\,, \tag{6.4}$$

*for any $\eta \in U$ if and only if for any $\eta \in U$, the complex curve $\psi_\eta : t + is \longmapsto s\Phi_t^\nabla(\eta)$, defined in a neighborhood of $0 \in \mathbb{C}$, is $J$-holomorphic.*

**Proof.** We observe first that the differential of the maps $\psi_\eta$ is given by

$$d_{t_0+is_0}\,\psi_\eta\Big(a\frac{\partial}{\partial t}+b\frac{\partial}{\partial s}\Big) \;=\; a\,d(s_0\,\mathbb{I}_{T_M})\dot{\Phi}_{t_0}^\nabla(\eta) + b\,T_{s_0\Phi_{t_0}^\nabla(\eta)}\,\Phi_{t_0}^\nabla(\eta)\,.$$

But

$$\begin{aligned}\dot{\Phi}_{t_0}^\nabla(\eta) &= \zeta^\nabla \circ \Phi_{t_0}^\nabla(\eta)\\ &= H^\nabla_{\Phi_{t_0}^\nabla(\eta)} \cdot \Phi_{t_0}^\nabla(\eta)\,,\end{aligned}$$

thanks to (6.1). Then using the property, $H^\nabla_{\lambda\eta} = d_\eta(\lambda\,\mathbb{I}_{T_M}) \cdot H^\nabla_\eta$, of the linear connection $\nabla$ (see the identity 8.5 in [Pali])we infer

$$d_{t_0+is_0}\,\psi_\eta\Big(a\frac{\partial}{\partial t}+b\frac{\partial}{\partial s}\Big) = \big(a\,H_{s_0\Phi_{t_0}^\nabla(\eta)} + b\,T_{s_0\Phi_{t_0}^\nabla(\eta)}\big)\cdot\Phi_{t_0}^\nabla(\eta)\,. \tag{6.5}$$

The complex curve $\psi_\eta$ is $J$-holomorphic if and only if

$$d_{t_0+is_0}\,\psi_\eta\Big(-b\frac{\partial}{\partial t}+a\frac{\partial}{\partial s}\Big) \;=\; J d_{t_0+is_0}\,\psi_\eta\Big(a\frac{\partial}{\partial t}+b\frac{\partial}{\partial s}\Big),$$

thus, if and only if

$$\big(-b\,H_{s_0\Phi_{t_0}^\nabla(\eta)} + a\,T_{s_0\Phi_{t_0}^\nabla(\eta)}\big)\cdot\Phi_{t_0}^\nabla(\eta) \;=\; J\big(a\,H_{s_0\Phi_{t_0}^\nabla(\eta)} + b\,T_{s_0\Phi_{t_0}^\nabla(\eta)}\big)\cdot\Phi_{t_0}^\nabla(\eta)\,.$$

For $s_0 \neq 0$ this is equivalent to (6.4). For $s_0 = 0$ this is equivalent to (6.3). We deduce the required conclusion. □

The condition (6.3) implies that $J$ is an $M$-totally real complex structure. We provide now the proof of corollary 1.3.

**Proof.** If we write $A = \alpha + iTB$ and $\alpha = H^\nabla - T\Gamma$, then $S := T^{-1}(H^\nabla - \overline{A}\,) = \Gamma + i\,B$. We set $S_k = \Gamma_k + i\,B_k$. From the proof of corollary 6.2 we know that in the case $J$ is integrable over $U$, the curve $\psi_\eta$ is $J$-holomorphic if and only if hold (6.4). The later rewrites as

$$H^\nabla_\eta \cdot \eta \;=\; -J_\eta T_\eta \cdot \eta.$$

Using the property

$$J_{\eta|\mathrm{Ker}\,d_\eta\pi} = -\alpha_\eta\,B_\eta^{-1}\,T_\eta^{-1},$$

(see the identity 1.7 in [Pali]), we infer that the previous identity is equivalent to

$$H^\nabla_\eta \cdot \eta = \alpha_\eta\,B_\eta^{-1}\cdot\eta. \tag{6.6}$$

Taking $d_\eta\,\pi$ on both sides of (6.6) we deduce $\eta = B_\eta^{-1}\cdot\eta$. Therefore (6.6) is equivalent to the system

$$\begin{cases} B_\eta \cdot \eta \;=\; \eta\,,\\[2mm] H^\nabla_\eta \cdot \eta \;=\; \alpha_\eta \cdot \eta\;. \end{cases} \tag{6.7}$$

Then the system (6.7) rewrites as

$$\begin{cases} \sum_{k\geqslant 1} B_k(\eta^{k+1}) \;=\; 0\,, \\ \\ \sum_{k\geqslant 1} \Gamma_k(\eta^{k+1}) \;=\; 0\;. \end{cases}$$

and thus as $S_k(\eta^{k+1}) \;=\; 0$ for all $k \geqslant 1$. We remind now that, according to the main theorem in [Pali], the integrability of the structure $J$ implies the condition $S_1 \in C^\infty(M, S^2 T_M^* \otimes_{\mathbb{R}} \mathbb{C}T_M)$. We infer $S_1 = 0$.

We notice that, with the notations of the statement of the main theorem in [Pali] the equation $\operatorname{Circ} \beta_k = 0$ hold for all $k \geqslant 1$. This combined with the identity $[\operatorname{Circ}, \operatorname{Sym}_{2,\ldots,k+2}] = 0$, implies

$$\operatorname{Circ} \operatorname{Sym}_{2,\ldots,k+2} \beta_k = 0\,, \tag{6.8}$$

for all $k \geqslant 1$. So if we apply the Circ operator to both sides of the definition of $S_2$ in the statement of the main theorem in [Pali], we infer $\operatorname{Circ} S_2 = \operatorname{Circ} \sigma_2 = 3\sigma_2$. If we evaluate this equality to $\eta^3$ we infer $S_2(\eta^3) = \sigma_2(\eta^3)$, which implies $\sigma_2 = 0$. We show now by induction that $\sigma_k = 0$ for all $k \geqslant 2$. Indeed by the inductive assumption

$$S_{k+1} \;=\; \frac{i}{(k+2)!} \operatorname{Sym}_{2,\ldots,k+2} \beta_k + \sigma_{k+1}\,.$$

Applying the Circ operator to both sides of this identity and using the equation (6.8), we infer $\operatorname{Circ} S_{k+1} = \operatorname{Circ} \sigma_{k+1} = 3\sigma_{k+1}$, which evaluated at $\eta^{k+2}$ gives $S_{k+1}(\eta^{k+2}) = \sigma_{k+1}(\eta^{k+2})$. We deduce $\sigma_{k+1} = 0$. Using the identity

$$\operatorname{Sym}_{2,\ldots,k+1} \operatorname{Sym}_{3,\ldots,k+1} = (k-1)!\, \operatorname{Sym}_{2,\ldots,k+1}\,, \tag{6.9}$$

we infer from the statement of the main theorem in [Pali] and with the notations there

$$S_k \;=\; \frac{i}{(k+1)!\,k!} \operatorname{Sym}_{2,\ldots,k+1} \theta_{k-1}\,,$$

for $k \geqslant 2$, with $\theta_1 := 2R^\nabla$ and

$$\begin{aligned} \theta_k \;:=\; & -2i\,(i d_1^\nabla)^{k-2} (\nabla R^\nabla)_2 \\ & + \sum_{r=3}^{k} (r+1)! \sum_{p=2}^{r-1} (i d_1^\nabla)^{k-r} (p S_p \wedge_1 S_{r-p+1})\,, \end{aligned}$$

for all $k \geqslant 2$. Moreover the equation $\operatorname{Circ} \beta_k = 0$, $k \geqslant 3$ rewrites as $\operatorname{Circ} \operatorname{Sym}_{3,\ldots,k+2} \theta_k = 0$. If we set $\Theta_k^\nabla := \theta_{k-1}$, for all $k \geqslant 2$ we obtain the required expansion.

On the other hand if the expansion in the statement of the lemma under consideration hold then $J$ is integrable thanks to the main theorem in [Pali] and $\operatorname{Circ} S_k = 0$ for all $k \geqslant 2$ ($S_1 = 0$). Indeed for $k = 2, 3$ this equality follows from the identities $\operatorname{Circ} \Theta_k^\nabla = 0$ and

$$[\operatorname{Circ}, \operatorname{Sym}_{2,\ldots,k+1}] = 0\,. \tag{6.10}$$

For $k \geqslant 4$ we use the identities (6.10), (6.9) and the equations satisfied by the covariant derivative operator $\nabla$. We deduce $S_k(\eta^{k+1}) \;=\; 0$, for all $k \geqslant 1$, which is equivalent to (6.6) and so to the fact that the curves $\psi_\eta$ are $J$-holomorphic. □